\documentclass[12pt]{amsart}
\usepackage{amsmath,amssymb}
\usepackage{graphicx,epsfig,subfigure,psfrag}
\begin{document}

%
%
\newtheorem{theorem}{Theorem}
\newtheorem{proposition}[theorem]{Proposition}
\newtheorem{lemma}[theorem]{Lemma}
\newtheorem{corollary}[theorem]{Corollary}
\newtheorem{definition}[theorem]{Definition}
\newtheorem{remark}[theorem]{Remark}
\numberwithin{equation}{section}
\numberwithin{theorem}{section}
\newcommand{\be}{\begin{equation}}
\newcommand{\ee}{\end{equation}}
\newcommand{\re}{\mathbb{R}}
\newcommand{\n}{\nabla}
\newcommand{\ren}{\mathbb{R}^N}
\newcommand{\iy}{\infty}
\newcommand{\pa}{\partial}
\newcommand{\ms}{\medskip\vskip-.1cm}
\newcommand{\mpb}{\medskip}
\newcommand{\BB}{{\bf B}}
\newcommand{\Am}{{\bf A}_{2m}}
\newcommand{\bL}{\BB^*}
\newcommand{\bLs}{\BB}
\renewcommand{\a}{\alpha}
\renewcommand{\b}{\beta}
\newcommand{\g}{\gamma}
\newcommand{\ka}{\kappa}
\newcommand{\G}{\Gamma}
\renewcommand{\d}{\delta}
\newcommand{\D}{\Delta}
\newcommand{\e}{\varepsilon}
\renewcommand{\l}{\lambda}
\renewcommand{\o}{\omega}
\renewcommand{\O}{\Omega}
\newcommand{\s}{\sigma}
\renewcommand{\t}{\tau}
\renewcommand{\th}{\theta}
\newcommand{\z}{\zeta}
\newcommand{\wx}{\widetilde x}
\newcommand{\wt}{\widetilde t}
\newcommand{\noi}{\noindent}
\newcommand{\lb}{\left (}
\newcommand{\rb}{\right )}
\newcommand{\lsb}{\left [}
\newcommand{\rsb}{\right ]}
\newcommand{\lab}{\left \langle}
\newcommand{\rab}{\right \rangle }
\newcommand{\gap}{\vskip .5cm}
\newcommand{\bz}{\bar{z}}
\newcommand{\bg}{\bar{g}}
\newcommand{\Ba}{\bar{a}}
\newcommand{\bt}{\bar{\th}}
\def\com#1{\fbox{\parbox{6in}{\texttt{#1}}}}

\title{\bf Very singular solutions \\ for thin film equations with absorption}

\author
{V.A. Galaktionov}

\address{
 Department of Math. Sci., University of Bath,
 Bath, BA2 7AY}
\email{vag@maths.bath.ac.uk}


\date{\today}



\begin{abstract}

The large-time behaviour of weak nonnegative and sign changing
solutions of the thin film equation (TFE) with absorption
\[
u_t = -\nabla \cdot (|u|^n \nabla \Delta u) - |u|^{p-1}u, 
\]
where $n \in (0,3)$
 and the absorption exponent $p$ belongs to the {\em subcritical} range
  $$
  \mbox{$
  p \in (n+1,p_0), \quad \mbox{with} \,\,\, p_0= 1+n + \frac 4 n,
   $}
   $$
   is studied. Firstly,
the standard free-boundary problem with zero-height, zero contact
angle and zero-flux conditions at the interface
   and bounded
compactly supported initial data is considered.  Very singular
similarity solutions (VSSs)  have  the form
 $$
 \textstyle{u_s(x,t) = t^{-\frac 1{p-1}} f(y), \quad
y= x/t ^\b, \quad \b = \frac{p-(n+1)}{4(p-1)}>0.}
 $$
  Here $f$
solves the quasilinear degenerate  elliptic equation
 $$
 \textstyle{- \n \cdot(|f|^n \n \D f) + \b \n f \cdot
y
+ \frac 1{p-1} \, f -|f|^{p-1}f=0} 
 $$
 that becomes an ODE for $N=1$ or in the radial setting in $\ren$.
By a combination of analytical, asymptotic, and numerical methods,
existence of various branches of similarity profiles $f$
parameterized by $p$ is established. Secondly,  in parallel,
changing sign VSSs of the Cauchy problem
 are
described.

This study 
is motivated by the detailed VSS results for
  the second-order porous medium equation  
with critical absorption ($u \ge 0$)
\[
\textstyle{ u_t = \n \cdot (u^n \n u) - u^p  \quad \mbox{in}
\,\,\, \ren \times \re_+, \quad p= 1+n + \frac 2 N,} \,\,\,\, n
\ge 0,
\]
which have been
  known since  the 1980s.

\end{abstract}

\maketitle

\section{\sc Introduction: very singular similarity solutions}
\label{Sect1}

\subsection{On the PME with critical
absorption}

Second-order quasilinear diffusion  equations with absorption are
typical for  combustion theory. The most well-known model, which
became a canonical object of intensive investigation in the
1970--90s, is the {\em porous medium equation} (PME) {\em with
absorption}
  \be
  \label{PME1}
 u_t = \n \cdot (u^n \n u) - u^p  \quad \mbox{in} \,\,\, \ren \times
 \re_+ \quad (u \ge 0),
   \ee
where $n \ge 0$ and $p \in \re$ are given parameters. Mathematical
theory of such PDEs was founded by  A.S.~Kalashnikov in the 1970;
see his fundamental survey \cite{Ka1} for the full PME history in
the 1950--80s. Besides new phenomena of localization and interface
propagation, for more than twenty years, the PME with absorption
(\ref{PME1}) was a crucial model for determining various
asymptotic patterns, which can occur for large times $t \gg 1$ or
close to finite-time extinction as $t \to T^-$ (for $p<1$). For
(\ref{PME1}), there are several
 parameter ranges with different asymptotics of solutions  such as
 \begin{gather*}
p>p_0=1+n + \textstyle{\frac 2N}, \quad p=p_0, \quad1+n < p <p_0,
\quad p=1+n, \\ 1< p < 1+n, \quad p=1, \quad 1-n < p<1, \quad p=
1-n, \quad p < 1-n,
\end{gather*}
 etc.; see details, references, and history in the book  \cite[Ch.~5, 6]{AMGV}.
A {\em transitional} behaviour for (\ref{PME1}) occurs precisely
at the first {\em critical absorption exponent}
 \be
\label{cr11} \textstyle{p_0= 1+n + \frac 2 N,}
 \ee
where logarithmically perturbed asymptotic patterns occur. VSSs,
as special key solutions,  are known to appear in the subcritical
range $p \in (n+1,p_0)$, while for $p>p_0$, the asymptotics of
solutions as $t \to \infty$ correspond to pure diffusive PME,
where the absorption term becomes negligible.

We will consider a new quasilinear parabolic model obtained by
adding to the standard thin film operator an extra absorption
term. This creates a non-conservative evolution PDE with the
differential operator that is only ``partially" divergent. We are
going to show that, nevertheless, this higher-order nonlinear PDE
exhibits several evolution features that are surprisingly similar
to those for the PME with absorption (\ref{PME1}). Moreover, we
also claim that mathematical features and properties of both
equations are also similar but, indeed, the results  are much more
difficult to justify for the higher-order case. Several
conclusions of the present paper remain interesting and  difficult open
mathematical problems, which we would like to stress upon.

We claim that many of the presented results and ideas cannot and will not be proved rigorously
in a few years to come,
but a reasonably detailed qualitative understanding of such popular nowadays models and nonlinear PDEs is essential.


\subsection{The TFE with absorption} Out goal is to
show that VSS phenomena exist in the new model that is the {\em
thin film equation} (TFE) {\em with absorption}
  \be
 \label{GPP}
 u_t = -\nabla \cdot (|u|^n \nabla \Delta u) - |u|^{p-1}u, 
   \ee
 where $n>0$ and $p>1$ are fixed exponents. For convenience,
it is written for solutions of changing sign to be also studied.
  Extra  absorption  terms in such TFEs  
are to model the effects of evaporation or permeability of the
surface. We refer to \cite{Oron97} and \cite{Beck05,Govor05} for
derivation of various TFE models including the non-conservative
cases. Other related details and references on such PDEs can be
found in the previous paper \cite{PetI} on the TFE (\ref{GPP})
devoted to the delicate case of the {\em critical} absorption
exponent
 \be
\label{cr222} \textstyle{p_0= 1+n + \frac 4 N.}
 \ee

In the present paper
 we study equation (\ref{GPP}) in the complementary
  {\em subcritical range}
 \be
 \label{scr1}
 \textstyle{
 n+1 < p < p_0= 1+n + \frac 4N.
  }
 \ee
 All necessary key references on various results on modern TFE theory
that are necessary for justifying principal regularity and other
assumptions on solutions will be presented below and are partially
available in \cite{PetI}. See also \cite[Ch.~3]{GSVR}, where
further references are given and several properties of TFEs with
absorption are discussed.

\smallskip

We will consider two main problems for the TFE (\ref{GPP}):

{\bf (I)} the standard free-boundary problem (FBP) admitting
nonnegative solutions, and

{\bf (II)} the Cauchy problem with another functional setting and
solutions of changing sign.

\smallskip

{\bf (I)} For the {\bf FBP}, equation (\ref{GPP}) is equipped
with {\em zero-height}, {\em zero contact angle}, and {\em
zero-flux} ({conservation of mass}) free-boundary conditions
 \be
 \label{GPP1}
 u=\n u=  {\bf \nu} \cdot (u^n \n \D u)  =0
 \ee
 at the singularity surface (interface) $\Gamma_0[u]$, which is
the lateral boundary of $ {\rm supp} \,u$ with the outward unit
normal ${\bf \nu}$.  Bounded, smooth, and compactly supported
initial data
  \be
 \label{u00}
 u(x,0)=u_0(x) \quad \mbox{in} \,\,\,\Gamma_0[u]\cap\{t=0\}
  \ee
 are added to complete a suitable functional setting of the
 FBP.
 We then assume that these three free-boundary conditions
  give a correctly specified problem for the fourth-order
parabolic equation, at least for sufficiently smooth and
bell-shaped initial data, e.g., in the radial setting; see
references below.

\smallskip

 {\bf (II)} For the {\bf Cauchy problem} (CP), we need solutions exhibiting the {\em
maximal} regularity at interfaces, and this demands oscillatory
character of such solutions. We refer to
 the book \cite[Ch.~3]{GSVR} and
\cite{Gl4, GBl6} for details concerning definitions and details on
oscillatory changing sign solutions of various TFEs and other
nonlinear higher-order PDEs.

\subsection{On the pure TFE and similarity solutions}

In what follows we will need well-known similarity solutions of
the
   usual  divergent  TFE with mass conservation
  \be
 \label{TFE1}
 u_t= - \n \cdot (|u|^n \n \D u).
   \ee
Earlier references
  on derivation of such fourth-order TFE
and related models
   can be
found in \cite{Green78, Smyth88}, where first analysis of some
self-similar solutions  for $n=1$ was performed. Source-type  
similarity solutions of (\ref{TFE1}) for arbitrary $n$ were
studied in \cite{BPelW92} for $N=1$ and in \cite{BFer97} for the
equation in $\ren$. More information on similarity and other
solutions can be found in \cite{BerHQ00, BerHK00, Bow01}.
 TFEs admit
non-negative solutions constructed by ``singular" parabolic
approximations of the degenerate nonlinear coefficients. We refer
to the pioneering paper \cite{BF1} and to various extensions in
\cite{Ell96, EllS, LPugh, WitBerBer} and the references therein.
 See also
\cite{Gr95} for mathematics of solutions of the FBP and CP of
changing sign (such a class of solutions of the CP will be
considered later on).

In both the FBP and the CP, the source-type solutions of the TFE
(\ref{TFE1}) take the form
 \be
\label{uss1} \textstyle{u_s(x,t) = t^{-  \b N} F(y), \quad y=
x/t^\b , \quad \mbox{with exponent} \,\,\,\b = \frac 1{4+nN},}
  \ee
 where, in the FBP,  $F(y) \ge 0$ is a radially symmetric compactly supported solution
of the ODE \cite{BPelW92, BFer97}
 \be
\label{Od11} {\bf B}(F) \equiv -\n \cdot (F^n \n \D F) + \b \n F
\cdot y + \b N F=0.
  \ee
 In the case $n=1$, the similarity profile for the FBP
is given explicitly,
 \be
\label{f112} \textstyle{F(y) = c_0 (a^2-|y|^2)^2, \quad c_0= \frac
1{8(N+2)(N+4)}, \quad a>0,}
  \ee
 which was first constructed in  \cite{Smyth88}.
Most advanced results of asymptotic stability for the TFE
(\ref{TFE1}) were obtained for the similarity solutions with
profiles (\ref{f112}) for $n=1$, where the rescaled equation is a
gradient system. See a full account of such studies in
\cite{CarrT02} ($L^1$-convergence), \cite{Car07}
($H^1$-convergence) and also  \cite{Gia08} for more general
properties including existence and uniqueness for the FBP (cf.
\cite[\S~6.2]{Gl4}, where these questions were treated in the von
Mises variables). For other values of $n \ne 1$, the results are
much weaker and are not complete. This emphasizes how difficult
the non fully divergent TFEs with non-monotone operators are for
the qualitative asymptotic study.

For the Cauchy problem, there exists a unique (up mass-scaling)
oscillatory similarity profile of (\ref{Od11}) for not that large
$n \in (0,n_{\rm h})$, $n_{\rm h}=1.759...$ \cite{Gl4}, while for
$n \in (0,1)$, existence and uniqueness are straightforward
consequences of the result in \cite{BMcL91} on oscillatory
solutions. See also some details in \cite[\S~3.7]{GSVR} and
interesting oscillatory features of similarity dipole solutions of
the TFE (\ref{TFE1}) observed in \cite{BW06}. Stability of such
sign changing similarity solutions is quite plausible but was not
proved rigorously being an open challenging problem.

 \subsection{On main results: comparison of logarithmically perturbed
 asymptotics for $p=p_0$ and VSSs for $p<p_0$}

 The critical case $p=p_0$ was studied in
\cite{PetI}, where it was shown that the asymptotic behaviour of
solutions is governed by a $\ln t$-perturbed similarity solution
(\ref{uss1}). For instance,
 for $n=1$, relying  on the explicit representation (\ref{f112})
and good spectral properties of the corresponding self-adjoint
linearised rescaled operator,
  it was shown that, for  $p=p_0= 2+
\frac 4 N$, the TFE with absorption (\ref{GPP}) admits asymptotic
patterns of the following form as $t \to \infty$:
 \be
\label{uln1} \textstyle{u(x,t) \sim (t \ln t)^{-\b N} F_*(x/t^\b
(\ln t)^{-\b N/4}) \quad \left(\b = \frac 1{4+N}\right),}
  \ee
  where
$F_*$ is a fixed rescaled profile from the family (\ref{f112})
with a uniquely chosen parameter $a=a_*(N)>0$.   For the
semilinear case $n=0$, i.e., for the fourth-order parabolic
equation
  written for solutions
 of changing sign
 \be
\label{sem1} u_t = - \D^2 u -|u|^{p-1}u,
   \ee
   the critical
behaviour like (\ref{uln1}) is known to occur at the
 critical exponent $p= 1+ \frac 4N$ \cite{GalCr},
  which is precisely (\ref{cr222}) with $n=0$.
   In this case, the centre manifold analysis also uses spectral
    properties of a non self-adjoint linear operator studied
     in \cite{Eg4}.   

\smallskip

 It turns out that, in the subcritical range
(\ref{scr1}),
 the generic behaviour of compactly supported solutions needs
 a different  similarity description leading to the VSSs.
 In the next
Section \ref{SectSub}, we describe the VSSs in the subcritical
range $p \in (n+1,p_0)$ for the FBP.
 In the final Section
 \ref{SectVCP}, we discuss
VSS structures   for the Cauchy problem admitting maximal
regularity solutions of changing sign. In general,  we show the
branching of the VSS profiles from the similarity patterns $F$ of
the pure TFE.

\section{VSSs for the FBP in the subcritical range $p \in (n+1,p_0)$}
\label{SectSub}

In the range (\ref{scr1}), the TFE with absorption (\ref{GPP})
admits the self-similar {\em very singular  solution} (VSS) of the
standard form
 \be
\label{Sub.1}
  \textstyle{u_s(x,t) = t^{-\frac 1{p-1}} f(y), \quad
y= x/t ^\b, \quad \b = \frac{p-(n+1)}{4(p-1)}>0,}  \ee
 where $f$
solves the nonlinear elliptic equation
 \be
\label{Sub.2}
 \textstyle{- \n \cdot(|f|^n \n \D f) + \b \n f \cdot
y
+ \frac 1{p-1} \, f -|f|^{p-1}f=0}. 
  \ee
  For
$n=0$ and $p \in (1, 1+\frac 4N)$, the existence of a finite
number of oscillatory profiles $f$ in $\ren$ corresponding to the
Cauchy problem for (\ref{sem1}) was detected in \cite{GW2} by
using a $p$-bifurcation analysis.

For the current FBP, the solvability of  (\ref{Sub.2}) and
existence of a nonnegative compactly supported $f$ is unknown even
in radial geometry, and we are going to present some analytic and
numerical evidence concerning this. In 1D, (\ref{Sub.2}) is an
ODE,
  \be
 \label{OD66}
 \textstyle{
 -(|f|^n f''')'+ \b f'y+ \frac 1{p-1}\, f -|f|^{p-1}f=0 \quad
 \mbox{for} \,\,\, y \in (0,y_0), \quad f'(0)=f'''(0)=0,
 }
   \ee
where we put two symmetry boundary conditions at the origin. At
the interface $y=y_0$, the solution expansion is as follows (see
\cite{BPelW92, Gl4, BFer97}):
  \be
\label{B112} \mbox{$
   f(y) = C_0(y_0-y)^2 +
 \frac{C_0^{1-n} \b  y_0}{(3-2n)(4-2n)(5-2n)} \,
   (y_0-y)^{5-2n}(1+o(1)),
   $}
  \ee
where $y_0>0$ and $C_0>0$ are arbitrary two parameters. Therefore,
in general, matching the 2D bundle (\ref{B112}), comprising two
positive parameters $y_0$ and $C_0$, with two symmetry conditions
in (\ref{OD66}) cannot give more than a countable number of
similarity FBP profiles $\{f_l\}$ (provided the parameter
dependence is analytic).

\begin{figure}
\centering
\subfigure[$p=2,$ $n=0.4$, $y_0\approx6.606$]{
\includegraphics[scale=0.4]{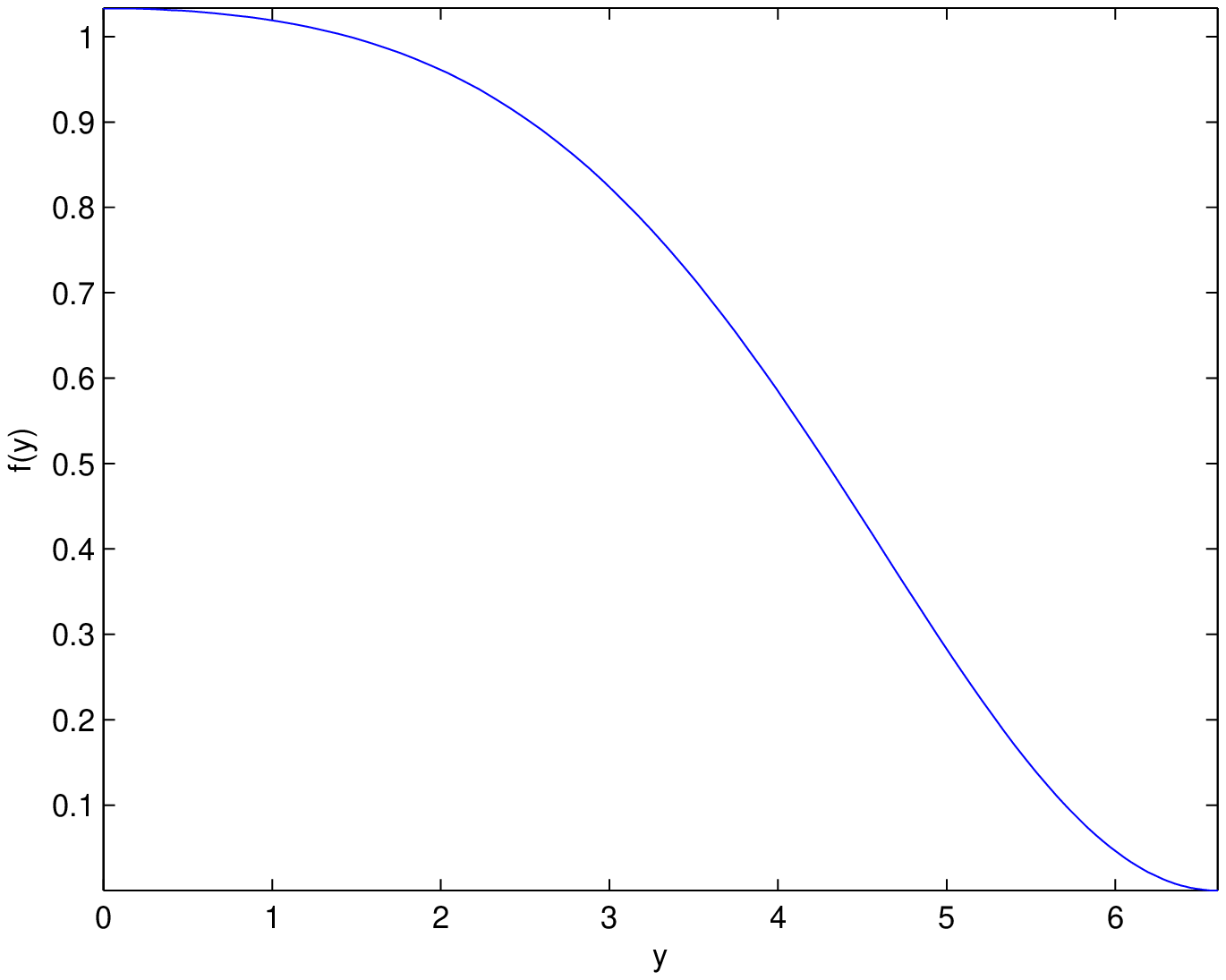}
\label{Fig6a}
}
\subfigure[$p=2,$ $n=0$, $y_0\approx5.109$]{
\includegraphics[scale=0.4]{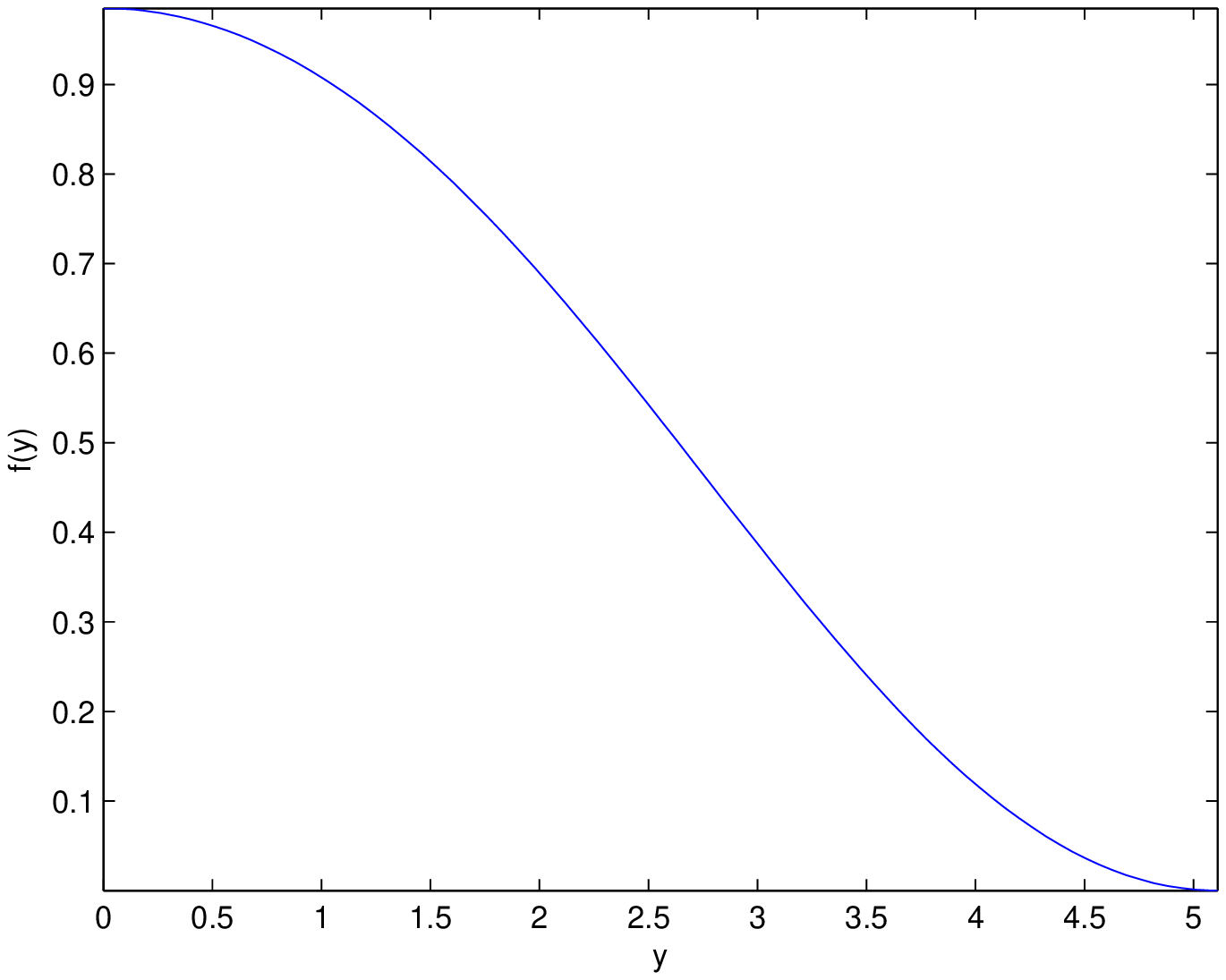}
\label{Fig6b}
}
\subfigure[$p=2,$ $n=0.8$, $y_0\approx14.822$]{
\includegraphics[scale=0.4]{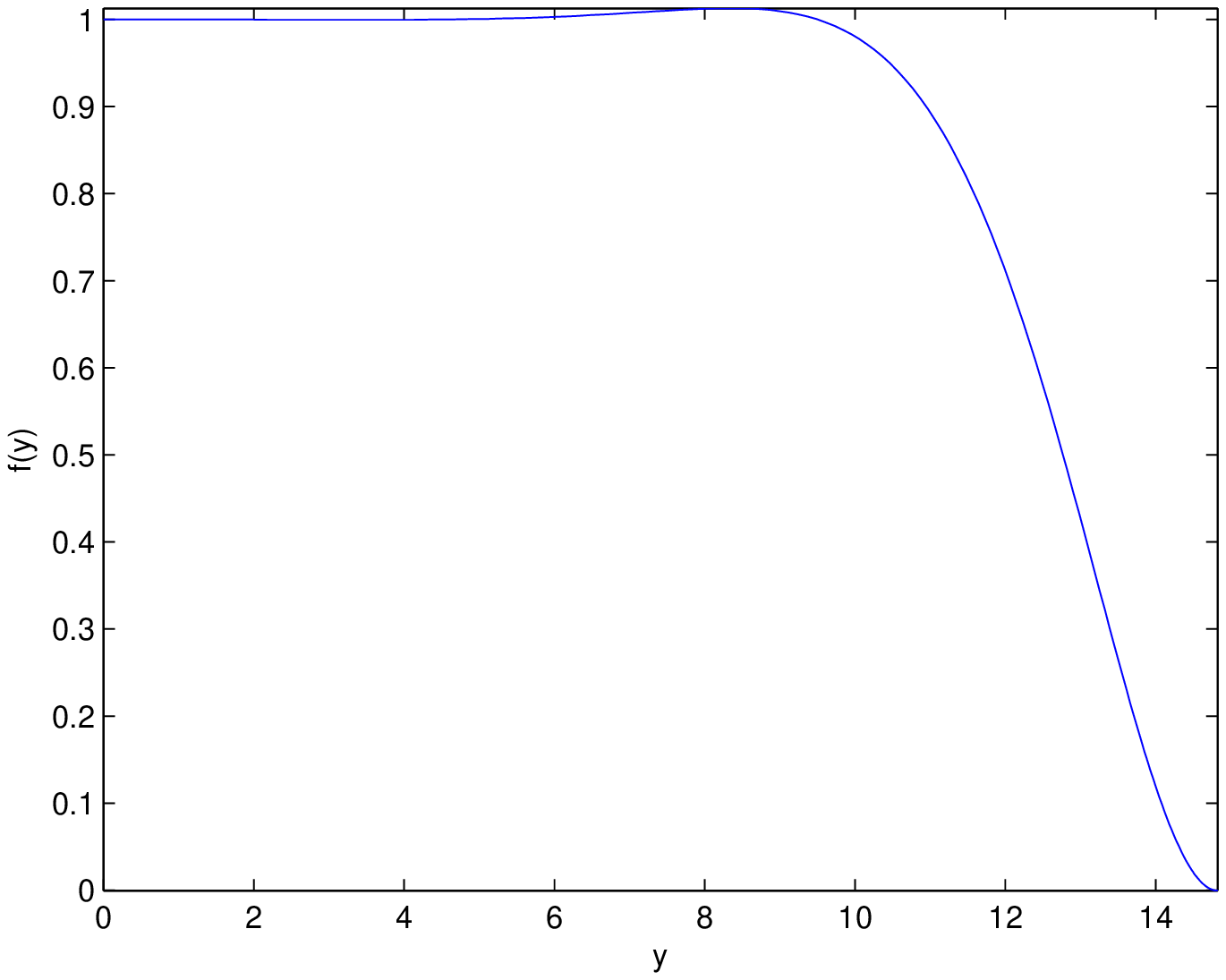}
\label{Fig6c}
}
\subfigure[$p=3,$ $n=1$, $y_0\approx 4.455$]{
\includegraphics[scale=0.4]{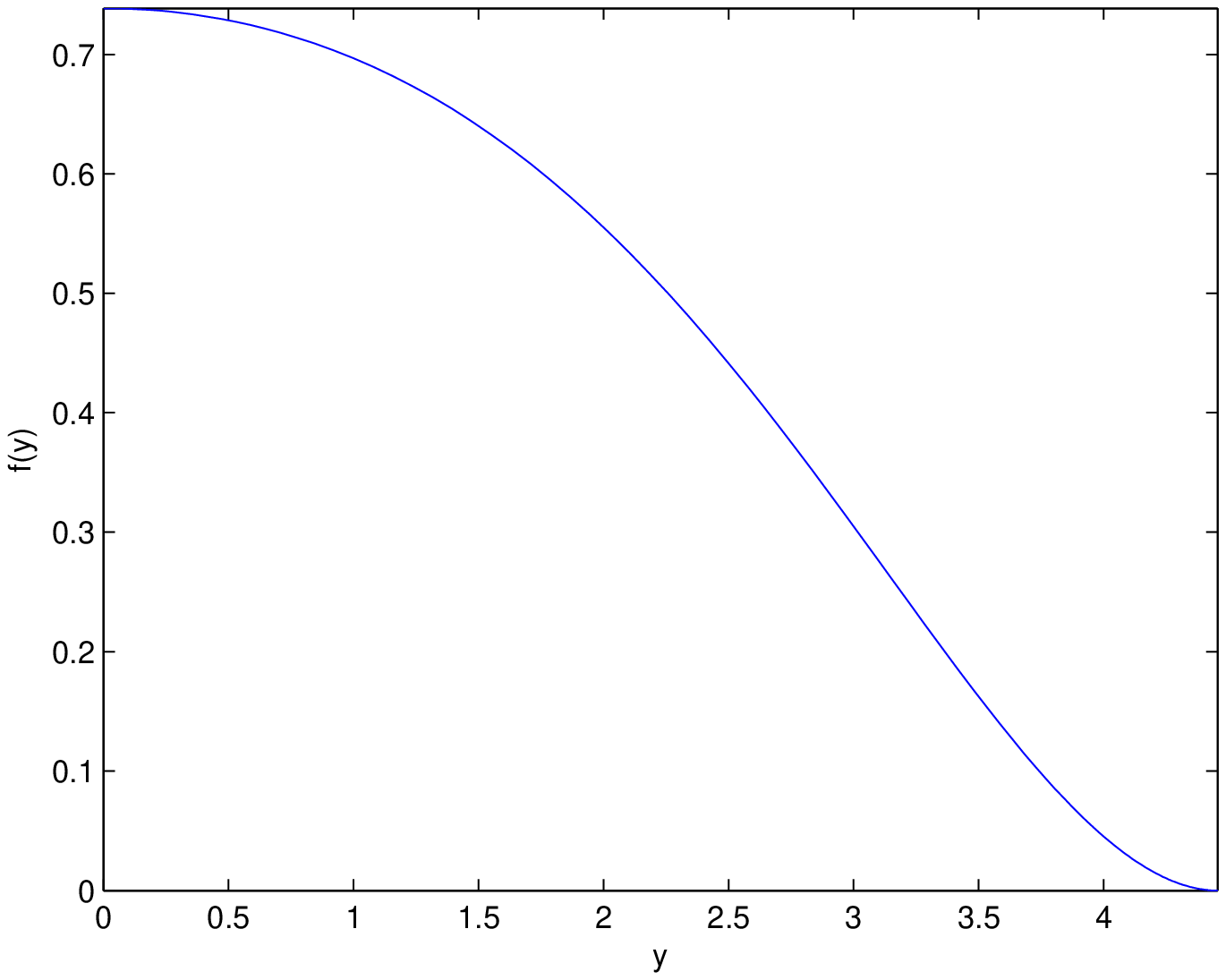}
\label{Fig6d}
}
\caption{FBP profiles satisfying
\eqref{OD66}, \eqref{B112} for various values of $p$ and $n$.}
\label{Fig6}
\end{figure}

Figure \ref{Fig6} shows  typical first positive symmetric VSS
profiles $f_0(y)$ constructed numerically by bvp4-solver in
MATLAB, by shooting in the $y_0$-parameter with conditions
  \be
 \label{sh22}
 f'(y_0)=f'''(y_0)=0.
   \ee
The correct choice of the interface location $y_0>0$ is obtained
from the zero-height condition $f(y_0)=0$, so in the limit, at
$y_0$, all three free-boundary conditions (\ref{GPP1}) are valid.
For numerics, we use the regularization in the fourth-order
operator in (\ref{OD66}) by replacing
 \be
 \label{Reg1}
  |f|^n \mapsto (\e^2 + f^2)^{\frac n2},
   \ee
   with typically $\e = 10^{-3}$, and similar tolerances of the
   method.
   Notice a big ``almost flat" part of the profile $f_0(y)$ in (c)
   for $p=2$ and $n=0.8$. Here, unlike other three cases, $p=2$ is
   sufficiently close for $n+1=1.8$, which is another critical
   exponent for the ODE (\ref{OD66}), since $\b=0$ then.

\begin{figure}
\centering \subfigure[$p=2, \, n=0.8$, $y_0=14.822...$]{
\includegraphics[scale=0.52]{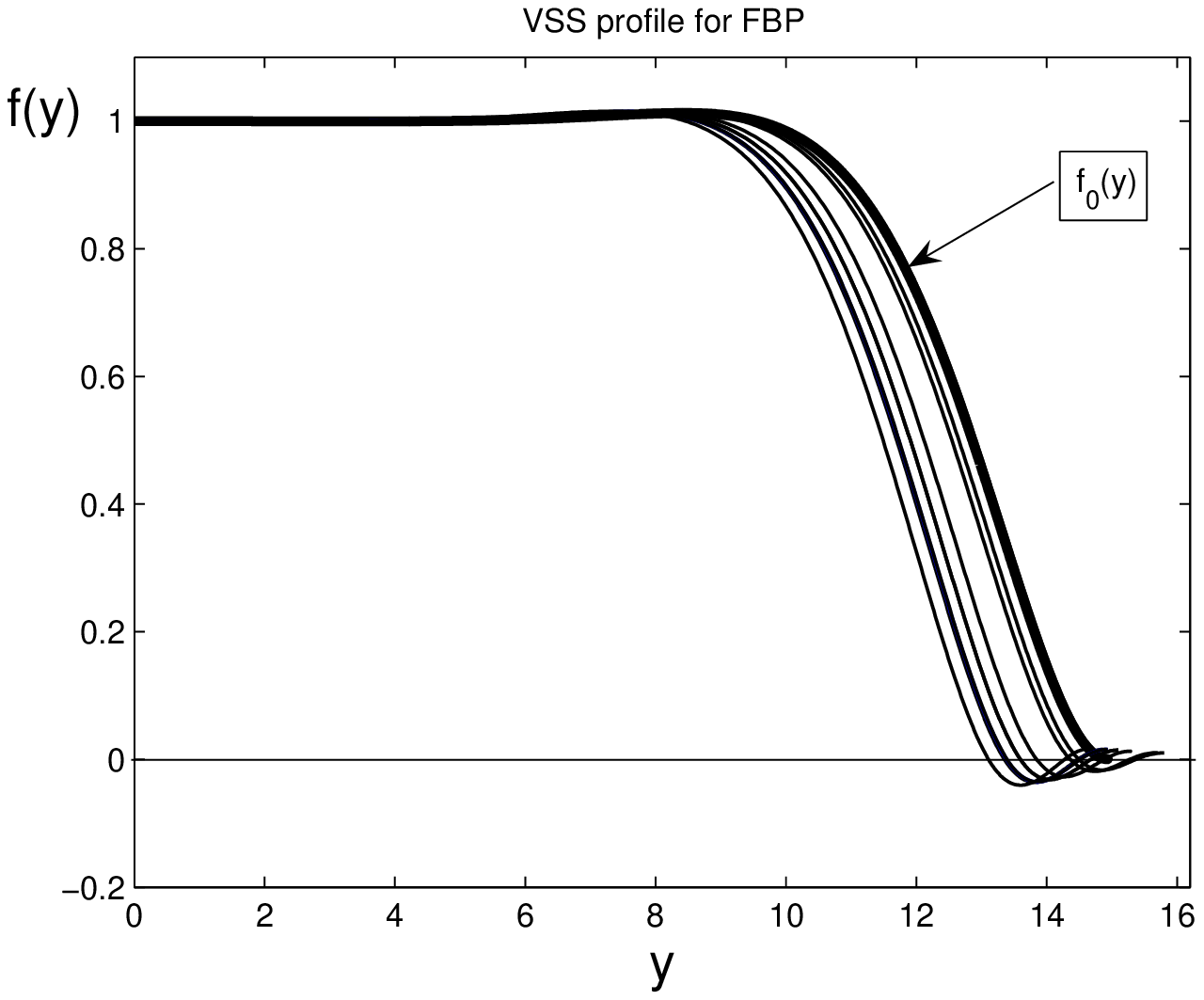}             
} \subfigure[$p=3, \, n=1$, $y_0=4.455...$]{
\includegraphics[scale=0.52]{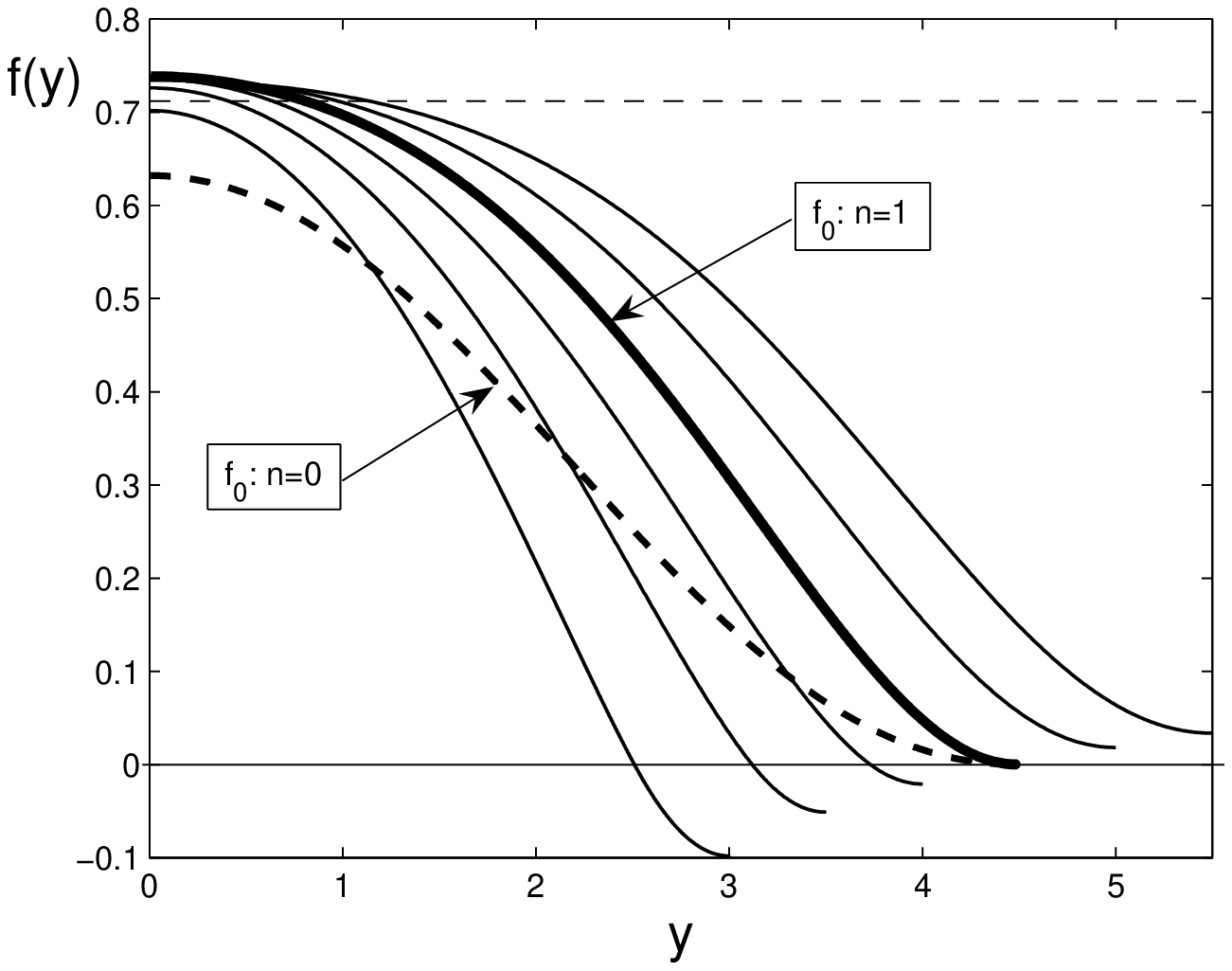}                        
}
 \vskip -.3cm
\caption{\rm\small Results of shooting of the profiles in Figure
\ref{Fig6} (c) and (d).}
 \label{Fig61}
\end{figure}

In Figure \ref{Fig61}, we explain further details and  the results
of the actual shooting procedure for the last two previous
patterns. For the reason of comparison, in (b) by dashed line we
show  the FBP profile for the semilinear case $n=0$ that was
studied in \cite{GW2} in the case of the Cauchy (not an FBP)
setting. It is curious that interfaces for $n=0$ and $n=1$ are
close to each other, but not the profiles.
 It is worth observing that Figure \ref{Fig61}(a)
reveals some oscillatory character of typical solutions near
interfaces. This indicates that there exist other patterns  $f_l$
for the FBP, which
 admit a few oscillations
near interfaces and form in the limit the solution of the Cauchy
problem; cf. \cite[Prop.~5.1]{PetI}.

Let us mention again that existence and multiplicity of VSS
profiles for (\ref{OD66}) are  open problems. In the next section,
we will discuss the question of existence of a finite number of
profiles in the subcritical range in the case of the CP.

 \vskip -.2cm





\section{VSSs in the Cauchy problem for $p \in (n+1,p_0)$}
 \label{SectVCP}

\subsection{On local oscillatory structure near interfaces}

The VSSs take the same self-similar form (\ref{Sub.1}), where the
radial rescaled profile $f$ of changing sign solves the ODE
(\ref{Sub.2}) in $\ren$.
 We refer to \cite[\S~5]{PetI} and \cite{Gl4} for details on the oscillatory
 structure of such similarity profiles close to finite interfaces,
 and also to \cite[Ch.~3]{GSVR}, where ``homotopy"-like approaches
 to the Cauchy problem for various TFEs are presented.

 Namely,
 it was shown
 that the asymptotic behaviour of $f(y)$ satisfying equation (\ref{OD66})
near the interface point as $y \to y_0^->0$ is given by the
expansion
 \be
\label{LC11}
 \mbox{$ f(y) = (y_0-y)^{\mu} \varphi(s) , \quad s
=\ln(y_0-y), \quad \mu = \frac 3 n,
 $}
   \ee
where, after scaling $\varphi \mapsto \b^{\frac 1 n} \varphi$, the
{\em oscillatory component}
 $\varphi$ 
 satisfies the following autonomous ODE, where exponentially small as $s \to -\infty$
  terms are omitted:
 \be
\label{m=2.11} 
 \textstyle{\varphi''' + 3(\mu-1) \varphi''  + (3
\mu^2 - 6 \mu +2) \varphi'+ \mu(\mu-1)(\mu-2)\varphi +  \frac
\varphi {|\varphi|^{n}}=0.}
  \ee
According to this singularity analysis, for $n \in (0,n_{\rm h})$,
where
 $$
   n_{\rm h}=1.759...
   $$
    is the {\em heteroclinic bifurcation point} for the ODE
    (\ref{m=2.11}),
there exists a stable (as $s \to +\infty$; at $s \to -\infty$ all
solutions are unstable in view of shifting in $y_0$) changing sign
periodic solution $\varphi(s)$ of (\ref{m=2.11}). According to
(\ref{LC11}), this gives similarity profiles of changing sign,
which being extended by $f(y) \equiv 0$ for $y>y_0$ forms a
 compactly supported solution  $f \in C^\a$ in a neighbourhood of $y=y_0$,  with $\a \sim \frac
 3n$. Notice that $\a \to +\infty$ as $n \to 0^+$, so the
 regularity at $y=y_0$ improves to $C^\infty$ at $n=0^+$.
 These functions  are
oscillatory near the interfaces. The first  results on the
oscillatory behaviour of similarity profiles for fourth-order ODEs
related to the source-type solutions of the divergent parabolic
PDE
  \be
  \label{mm1}
  u_t = - (|u|^{m-1} u)_{xxxx} \quad (m>1),
   \ee
 were obtained in \cite{BMcL91}.
 It turns out that these results can be applied to the rescaled TFE
 (\ref{Od11}), but for $n \in (0,1)$ only (the ODEs for
 (\ref{TFE1}) and (\ref{mm1}) then coincide after change).
 Some existence and multiplicity of
periodic solutions of ODEs such as (\ref{m=2.11}) are known
\cite{Gl4}, and often lead to a number of open problems.
Therefore,  numerical and some analytic evidence remain key,
especially for sixth and higher-order TFEs, \cite[Ch.~3]{GSVR}.

Thus, we express the above  results as follows: there exists a 2D
bundle of asymptotic solutions of (\ref{OD66}) in $\re$ having the
expansion (\ref{LC11}),
 \be
 \label{2par}
 f(y)=(y_0-y)^{\frac 3n} \varphi(s+s_0),
  \quad \mbox{with two parameters $y_0>0$ and $s_0 \in \re$},
  \ee
  where we also take into account the phase shift $s_0$ of the periodic orbit
 $\varphi(s)$.

Therefore, matching the 2D bundle (\ref{2par}) with two symmetry
conditions at the origin in (\ref{OD66}) leads to a reasonable
well-posed problem of 2D matching, which remains essentially open
still. In the case of analytic dependence on parameters involved,
such a problem cannot possess more than countable set of isolated
solutions. Actually, we are going to show that the number of VSS
profiles is always finite.

\subsection{Global behaviour of VSS profiles}

\begin{figure}
\centering \subfigure[$f_0(y)$ for $p=2$]{
\includegraphics[scale=0.52]{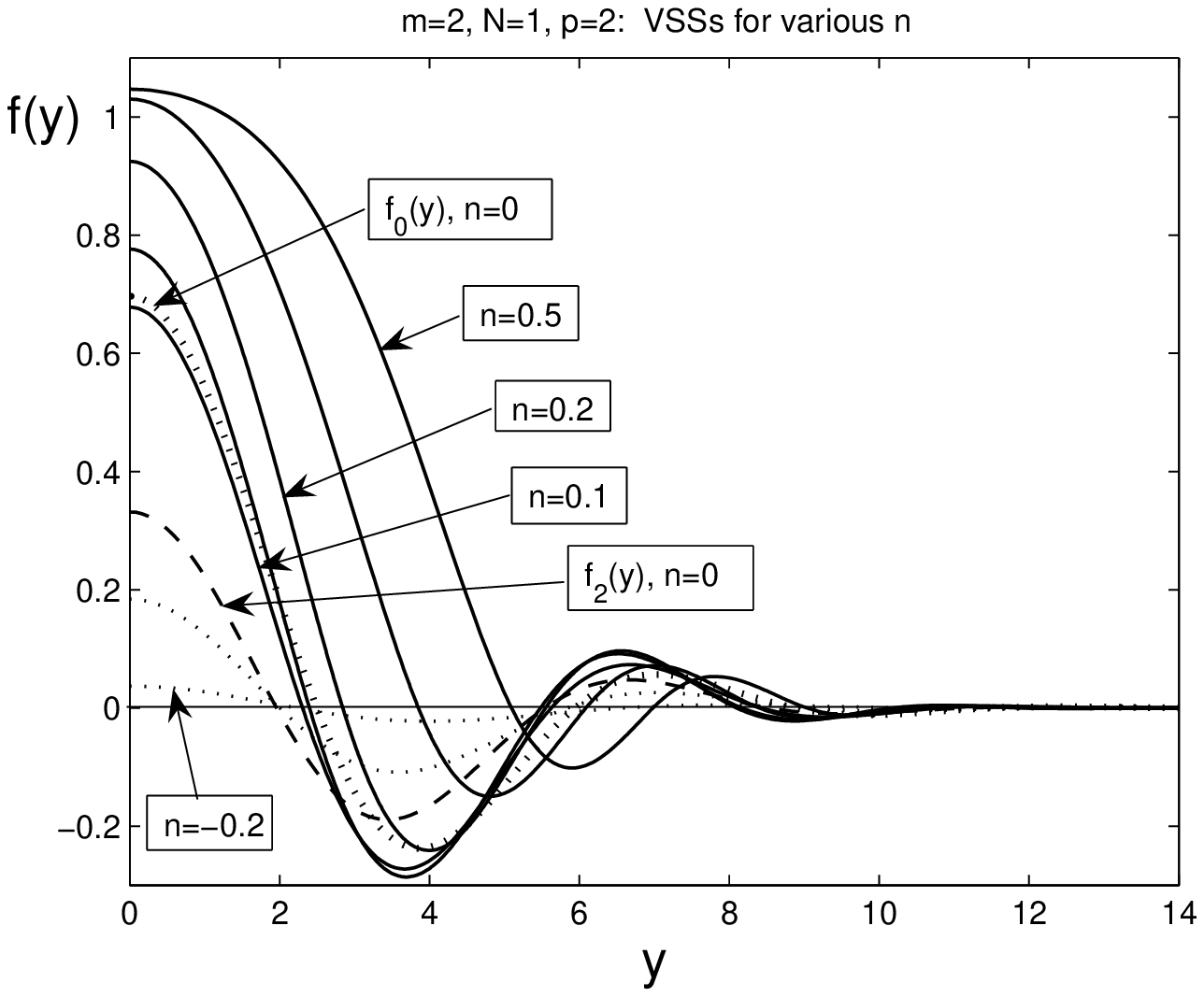}    
} \subfigure[$f_0(y)$ for $p=3$]{
\includegraphics[scale=0.52]{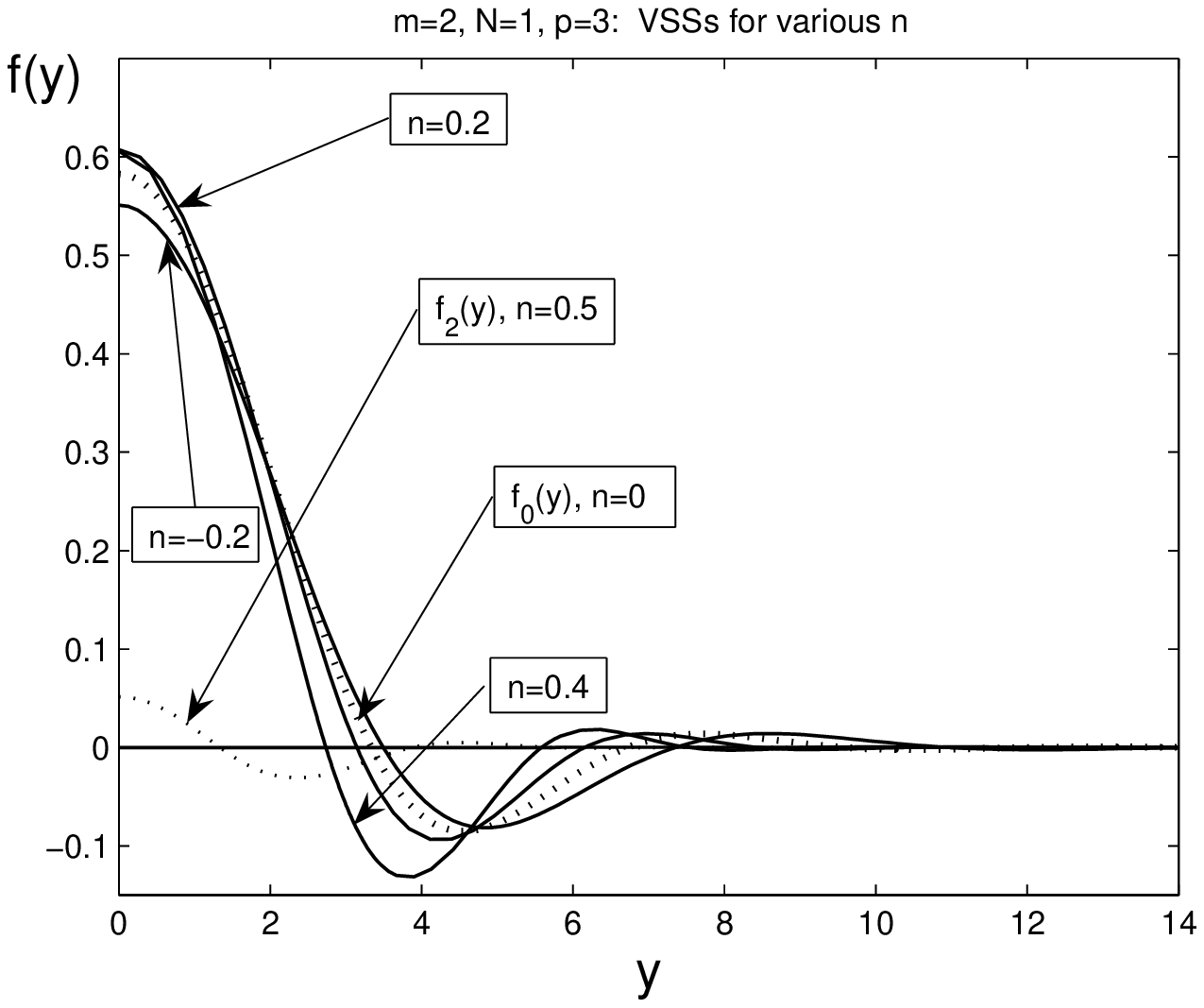}               
}
 \vskip -.4cm
\caption{\rm\small First oscillatory VSS profiles $f_0$ and $f_2$
of the CP satisfying (\ref{Sub.2}) in $\re$.}
 \label{F032}
\end{figure}

In Figure \ref{F032}, we present the first even VSS profile
$f_0(y)$ (and some $f_2(y)$) satisfying (\ref{Sub.2}) for $N=1$,
$p=2$, or 3 and various $n $.
 Here we fix the same regularization (\ref{Reg1}) where $\e$ and
 Tols are about $10^{-2}$, which is sufficient accuracy.
Notice that in 1D,
 $$
 \begin{matrix}
 f_0, \, f_2, \, f_4,... \,\, \mbox{are even functions, and} \smallskip \\
 f_1, \, f_3, \, f_5,... \,\, \mbox{are odd}. \qquad \qquad \quad
 \,\,
 \quad
  \end{matrix}
   $$

As we have mentioned,
 the case $n=0$ corresponds to smooth VSSs  for the semilinear
parabolic equation (\ref{sem1}),
 which were studied in \cite{GW2}, so we can always compare the results
  with those for the TFE (\ref{GPP}) with small $n>0$.
For $n=0$, the VSS profiles $\{f_l, \, l \ge 0\}$ are known to
appear at subcritical (for $p < p_l$) pitchfork bifurcations at
critical exponents
 \be
 \label{crN}
 \mbox{$
 p_l= 1+\frac 4{N+l}, \quad l=0,1,2,... \, .
  $}
  \ee

Figures  \ref{F032}(a) and (b) show a strong similarity of the
corresponding  VSS profiles $f_0$ for various $n$, and, moreover,
confirm that the
profiles can be continuously  
 deformed to each
other as $n \to 0^+$. This is related to a general homotopy
approach to the Cauchy problem for TFEs and other degenerate PDEs
with non-monotone operators; see \cite[Ch.~3]{GSVR} and \cite{Gl4,
GBl6}.

In Figure \ref{F032}(a) for $p=2$, we also show the smaller VSS
profile $f_2(y)$ for $n=0$ (the dashed line), and also a couple of
profiles for negative $n=-0.1$ and $-0.2$. These are $f_2(y)$, and
not $f_0(y)$. In (b), we also calculate $f_0$ for the negative
$n=-0.2$. Recall that for $n<0$, (\ref{GPP}) demonstrates typical
features of a {\em fast-diffusion} problem. There is no finite
propagation in this case, but solutions are equally oscillatory as
$y \to + \infty$, thus inheriting this property from $n=0$.
Notice that for $n=0$, according to (\ref{crN}),
 $$
  \mbox{$
 p_2= 1 + \frac 43 <3,
  $}
  $$
so that $f_2(y)$ does not exist for $p=3$. But $f_2$ exists for
positive $n=0.5$ (the dotted line in (b)) and is sufficiently
small, meaning that the corresponding critical bifurcation
exponent $p_2(n)$ is slightly above $3$.

\subsection{On a boundary layer as $p \to n+1$}

 For $n=0$ and any $p
\in (1,p_0)$, there exists a finite number
 $$
 \mbox{$
 M \sim \frac {p_0-p}{p-1} \to +\infty \quad \mbox{as} \,\,\, p
\to 1^+ \quad (n=0),
 $}
 $$
 of different VSS profiles, which are obtained by standard
 bifurcation theory, \cite{GW2}. We expect that a similar
 multiplicity property remains valid for the TFE for $n>0$, though
 in this case bifurcation branches  are
 governed by the linearised TFE operator
 (see \cite[\S~2]{PetI} and the results below),
   so that
 bifurcation points $\{p_l\}$ are not given explicitly by a discrete
 spectrum $\{\l_l\}$ of a non self-adjoint operator as in the
 semilinear case $n=0$; cf. \cite[Lemma~4.1]{GW2}.
We analyze this kind of bifurcation in the next subsection.

\begin{figure}
\centering
\includegraphics[scale=0.65]{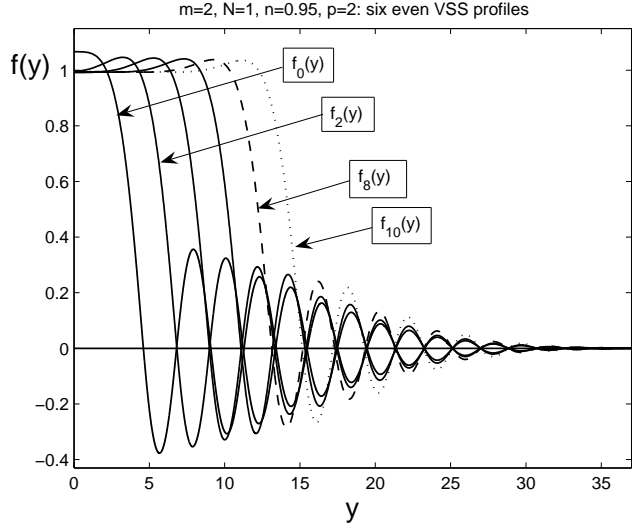}                     
 \vskip -.4cm
\caption{\rm\small Oscillatory behaviour of the VSS profiles for
$p \approx n+1$ increases: first six even VSS solutions for $p=2$
and $n=0.95$.}
 \label{Def1}
\end{figure}

A typical strong oscillatory behaviour  
of the VSS profiles for  $p \approx (n+1)^+$ is shown in Figure
\ref{Def1}, for $p=2$ and $n=0.95$. We present here first six even
VSS profiles from different $p_{2l}$-branches; see further
explanations below.
 Notice
 formation of
an interesting ``boundary layer" as $p \to n+1$, where the VSS
profiles $f(y)$ become more and more oscillatory reflecting the
fact that a suitable solution satisfying at $p=n+1$ the ODE
(\ref{OD66}) for $N=1$,
  \be
 \label{N1}
 \mbox{$
-(|f|^nf''')'+ \frac 1n \, f - |f|^n f=0 \quad \mbox{in} \,\,\,
\re, $}
  \ee
 does not exist. Numerics in Figure \ref{Def1} suggest  that solutions of (\ref{N1}) are highly
 oscillatory and are not compactly supported. Recall that the
 oscillatory structure (\ref{LC11}) near interfaces demands
 the linear term $+\b f'y$ in the ODE (\ref{OD66}) with $\b >0$.
 For $\b=0$,  such solutions do not exist.
Moreover, for $\b<0$, there exist positive solutions  near
interfaces, which correspond to blow-up problems  \cite{Bl4}.

\subsection{Bifurcation of the first $p$-branch at  $p = p_0^-$: a nonlinear version}

 We now study
the behaviour of the first $p$-branch of the VSS profiles
$f=f_0(y)$, when $p$ approaches from below the critical exponent
 (\ref{cr222}).
For $n=0$, such a behaviour was studied in \cite{GW2} by a
standard bifurcation approach showing that the VSS profiles vanish
with the rate:
  \be
\label{Dec1} \|f\|_\infty \sim (p_0-p)^{\frac N4} \to 0 \quad
\mbox{as} \,\,\, p \to p_0^-.
  \ee
The bifurcation analysis in \cite{GW2} was based on known point
spectrum (\ref{ss5}) and other spectral properties of a
non-self-adjoint linear operator ${\bf B}$, which is  (\ref{Od11})
for $n=0$. The  operator ${\bf B}$ appears in the ODE (\ref{Od11})
for $n=0$ generating the rescaled kernel of the fundamental
solution of the bi-harmonic operator $D_t+ \D^2$.

\smallskip

\noi{\bf Local bifurcation from $p_0$.} We now perform a formal
{\em nonlinear} version of such a bifurcation (branching) analysis
for $n>0$. As usual, according to classic branching theory
\cite{KrasZ, VaiTr}, a justification (if any) is performed for the
equivalent quasilinear integral equation with compact operators.
For simplicity, we present basic computations for the differential
version.

We introduce the small parameter $\e= p_0-p$, so that, as $\e \to
0$,
 $$
 \mbox{$
 \frac 1{p-1}= \frac N{4+nN} + \e \frac{N^2}{(4+nN)^2} + O(\e^2)
 \,\,\, \mbox{and} \,\,\,
  \b=\frac 1{4+nN} - \e \frac{nN^2}{4(4+nN)^2} + O(\e^2).
   $}
   $$
Substituting this expansion into (\ref{Sub.2}) and performing the
standard linearization yields
  \be
 \label{eq4}
  \begin{matrix}
 {\bf B}(f) + \e ({\mathcal L}_1 f - |f|^{n+\frac 4N} f \ln|f|)
 -|f|^{n+\frac 4N} f+ ... =0,\smallskip\smallskip\\
 \mbox{where} \quad  {\mathcal L}_1=
 \frac {N^2}{(4+nN)^2}(N I - \frac n4 \, y \cdot \nabla)
  \end{matrix}
  \ee
  is a linear operator, and ${\bf B}$ is the rescaled operator (\ref{Od11}) of the pure
   TFE.
Notice that the fact that the operator ${\bf B}$ in (\ref{eq4})
occurs in the rescaled
 pure TFE correctly describes the essence of a
 ``nonlinear bifurcation phenomenon" to be revealed.

Next, we use an extra invariant scaling of operator ${\bf B}$  by
setting
  \be
 \label{ep1}
 f(y)=b \tilde F(y/b^{n/4}),
  \ee
 where $b=b(\e)>0$ is a small parameter, $\b(\e) \to 0$ as $\e \to 0$, to be determined.
 Substituting (\ref{ep1}) into (\ref{eq4}) and omitting
 higher-order terms yields
  \be
 \label{ep2}
 {\bf B}(\tilde F) + \e {\mathcal L}_1 \tilde F- b^{n+\frac 4N}
|\tilde F|^{n+\frac 4N} \tilde F+...=0.
  \ee

 Finally, we perform linearization $\tilde F=F+Y$, where $F$ is
 the  ``fundamental",  supported in $B_1$, similarity profile
 of the Cauchy problem for the pure TFE, i.e., satisfying (\ref{Od11})
 for $N=1$. This yields the
 non-homogeneous problem
   \be
  \label{ep3}
  {\bf B}'(F)Y + \e {\mathcal L}_1 F - b^{n+\frac 4N}
|F|^{n+\frac 4N}  F+...=0.
  \ee
 Here the derivative is given by
  $$
   \mbox{$
  {\bf B}'(F)Y= - \n \cdot [|F|^n (\frac n F \, (\n \D F)Y + \n \D
  Y))] + \b y \cdot \n Y + \b N Y.
   $}
   $$

 The rest of the analysis depends on assumed good  spectral properties of
 the linearised operator ${\bf B}'(F)$. We follow the lines of a
 similar analysis performed for the FBP case in
 \cite[\S~2]{PetI}, where, in the FBP for $n=1$ and $F$ given
 by (\ref{f112}), the operator ${\bf B}'(F)$ turns out to possess   a
 (Friedrichs) self-adjoint  extension with compact resolvent and discrete
 spectrum. Such a self-adjoint extension does not exist for the oscillatory $F(y)$.
 Here we need to use general theory of non-self-adjoint operators;
  see e.g.,
 \cite{GGK}.
  A proper functional setting of this
 operator is more straightforward for $N=1$
 (and in the radial setting), where, using the behaviour of
 $F(y) \to 0$ as $ y \to 1$, it is possible to check whether the
 resolvent is compact is a suitable weighted $L^2$ space. In
 general, this is a difficult problem; see below.

We assume that such a proper functional setting is available
 for ${\bf B}$, so
we deal with operators having solutions with ``minimal"
singularities at the boundary $S_1$, where the operator is
degenerate and singular.
  Namely,  we find the first
 eigenfunction $\psi_0$ with $\l_0=0$ of ${\bf B}'(F)$.
 Let $\psi_0^*$ be the corresponding first eigenfunction
  of the adjoint operator $({\bf B}')^*(F)$ defined in  a natural
  way using the topology of the dual space $L^2(B_1)$ and having
  the same point spectrum  (the latter  is true for compact operators in a suitable space
  \cite[Ch.~4]{KolF}).
 Moreover,  it can be seen from the
 divergent form of the linearised operator ${\bf B}'(F)$ that, after bi-orthonormalisation,
 $$
  \langle \psi_\a, \psi_\b^* \rangle= \d_{\a \b},
  $$
  the
corresponding first eigenfunction of $({\bf B}')^*(F)$ can be
taken as
  \be
 \label{ep4}
 \psi_0^*(y) \equiv 1.
   \ee
   This simplifies the rest of computations, though one can restore these
   similarly for arbitrary $\psi_l^*$, as suggested later on for
   finding other critical bifurcation points $\{p_l\}$.

 Further,   we assume that there exists the orthogonal subspace ${\rm
   Span}\{\psi_l, l \ge 1\}\bot \psi_0$ of eigenfunctions of
   ${\bf B}'(F)$, and we look for solutions of (\ref{ep3}) in the form
    $$
    Y= C \psi_0+ w,
     $$
      with a constant $C=C(\e)$ and a function $w \bot \psi_0$, i.e.,
      $\langle w, \psi_0^*\rangle =0$. Recall again that, in doing so,
    we need to transform (\ref{ep3}) into an equivalent integral equation
    with compact operators, but for convenience, we continue our
    computations using the differential version; see some  details
   in \cite[\S~3]{GW2}.

    Thus,
multiplying (\ref{ep3}) by $\psi_0^*$ in
$L^2(B_1)$ and integrating by parts in the differential term $y
\cdot \n F$ in ${\mathcal L}_1 F$, we obtain the following
orthogonality condition
 of solvability 
(Lyapunov-Schmidt's branching equation \cite[\S~27]{VaiTr}):
  \be
 \label{ep41}
  \mbox{$
   \e\frac{N^2}{4(4+nN)} \int F= b^{n+\frac 4N} \int |F|^{n+\frac 4N}
   F.
   $}
    \ee
Therefore, the parameter $b(\e)$ in (\ref{ep1}) for $p \approx
p_0^-$ is given by (cf. (\ref{Dec1}) for $n=0$)
  \be
 \label{ep5}
 \mbox{$
 b(\e)= [\g_0 (p_0-p)]^{\frac N{4+nN}}, \quad \mbox{with}
 \,\,\,\, \g_0 = \frac{N^2}{4(4+nN)} \int F \big/\int |F|^{n+\frac 4N} F,
 $}
     \ee
    provided that $\int |F|^{n+\frac 4N} F>0$ (not an easy
    inequality that can be checked numerically).

\smallskip

    For $n=0$, a rigorous justification of this bifurcation
    analysis can be found in \cite[\S~6]{GW2}, where a
    countable number of $p$-branches originated at bifurcation
    points (\ref{crN})
     was detected
 on the basis of known spectral properties of the
 corresponding linear operator (\ref{ss4}); see details in
  \cite{Eg4}.
    For $n>0$, the justification needs spectral properties of the
    linearised operator ${\bf B}'(F)$ and the corresponding
    adjoint one $({\bf B}'(F))^*$, which
 remain an open problem.
 In particular, it would be important to
  know that the bi-orthonormal eigenfunction
  subset $\{\psi_l\}$ of the operator ${\bf B}'(F)$
  is complete and  closed in a weighted $L^2$-space (for $n=0$, such results are available \cite{Eg4}).   We expect that for $n \approx 0$, there
  exist critical exponents for the TFE with absorption that are
 close to those in (\ref{crN}) at $n=0$. This can be checked by
 standard
 branching-type calculus; see Appendix in \cite{GHCo},
 where nonlinear eigenfunctions of the rescaled PME in $\ren$ were
 studied by a branching approach.

\smallskip

\begin{figure}
\centering
\includegraphics[scale=0.7]{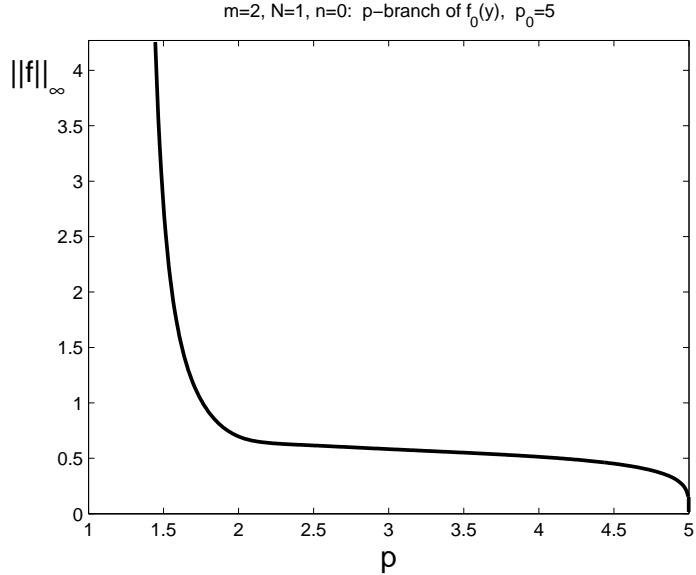}
 \vskip -.4cm
\caption{\rm\small $N=1$ and $n=0$: the first $p$-bifurcation
branch that is originated at the first critical exponent $p_0=5^-$
and blows up as $p \to 1^+$.}
 \label{Fn0}
\end{figure}

\noi{\bf On global extension of $p$-branches.}  For comparison, we
begin with Figure \ref{Fn0} that presents the first monotone
branch of VSS profiles $f_0(y)$ in the semilinear case $n=0$ that
exists for all $p \in (1,p_0=1+\frac 4N)$; cf. \cite{GW2}.
  On the
vertical axis, we put $\|f\|_\infty$ that, for such profiles, is
simply $f(0)$.

\begin{figure}
\centering \subfigure[$n=0.11$]{
\includegraphics[scale=0.52]{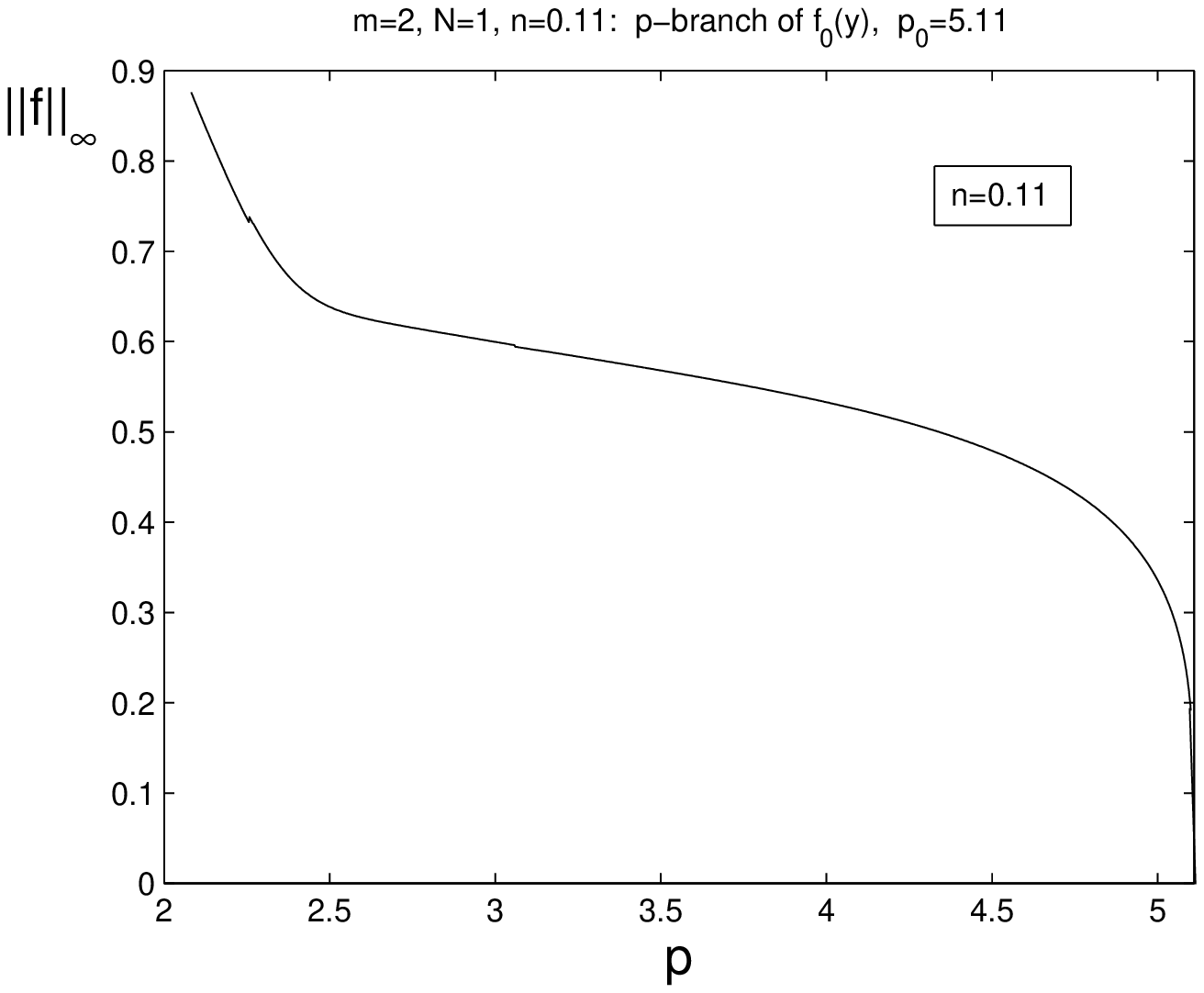}     
} \subfigure[$n=0.13$]{
\includegraphics[scale=0.52]{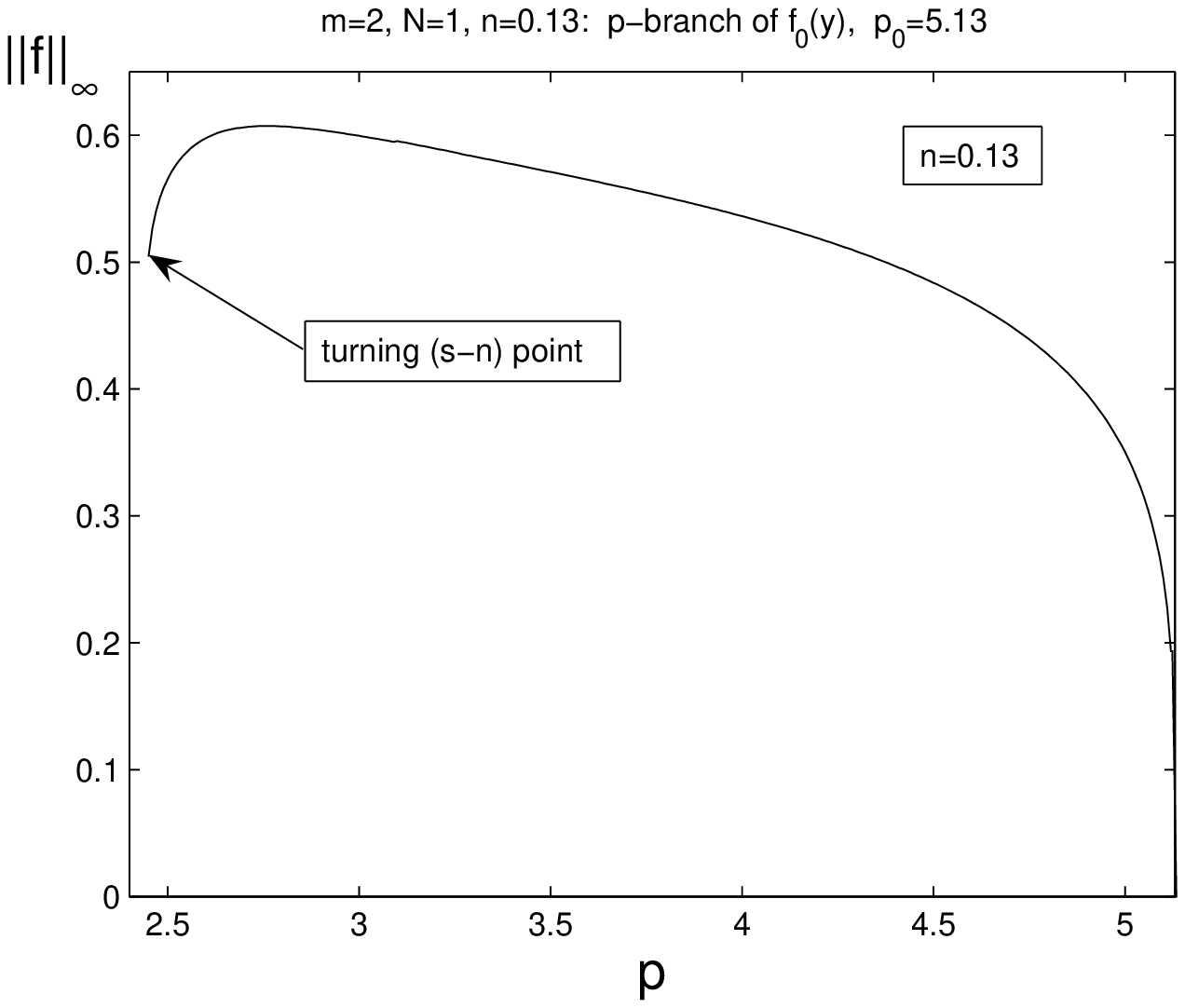} 
}
 \vskip -.4cm
\caption{\rm\small $p$-branches of VSS profiles for $N=1$:
$n=0.11$ (a) and $n=0.13$ (b).}
 \label{FF12}
\end{figure}

Figure \ref{FF12} show that such monotone branches persist until
$n=0.11$ (a), while, for a slightly larger $n=0.13$, we first
observe a non-monotone branch of patterns $f_0(y)$ (b), and this
persists for most of larger $n$'s. This is a {\em quasilinear}
phenomenon to be discussed in greater detail below.
 Moreover, we observe a typical {\em turning} point of the
$p$-branch at $p= 2.45...$, which, as usual,  is characterized by
existence of a non-trivial centre subspace of the linearized
operator,
 $$
 0 \in \s({\bf B}'(f)).
  $$
 These  global bifurcation diagrams
are calculated with the enhanced Tols = $10^{-4}$ and small step
sizes $\D p \sim 10^{-3}$. We claim that the turning, saddle-node
bifurcation for $N=1$ occurs above the critical exponent
 $$
 n_{\rm s-n} \approx 0.12.
  $$

\begin{figure}
\centering
\includegraphics[scale=0.7]{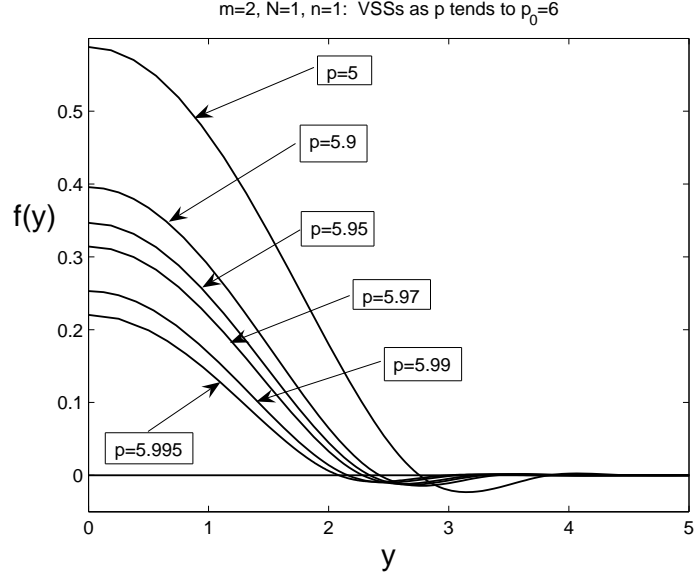}                 
 \vskip -.4cm
\caption{\rm\small $N=n=1$:  VSS profiles vanish   as $p \to p_0^-
=6$.}
 \label{FBif1}
\end{figure}

\begin{figure}
\centering
\includegraphics[scale=0.7]{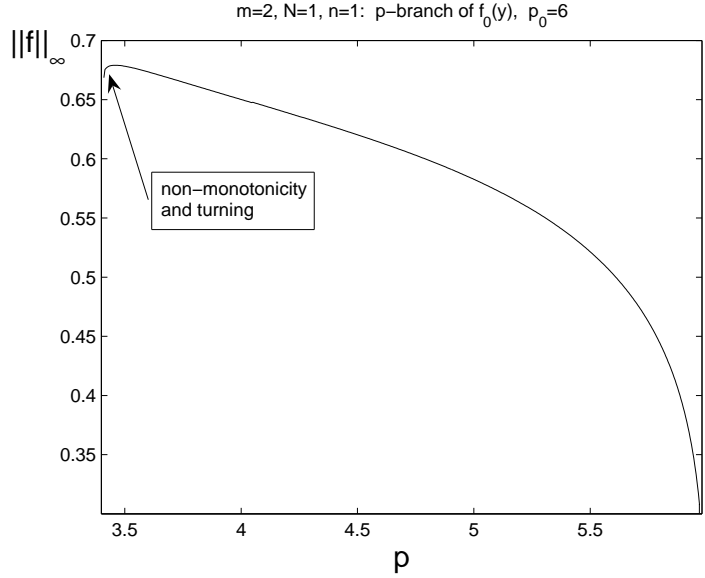}            
 \vskip -.4cm
\caption{\rm\small $N=n=1$: the first $p$-bifurcation branch of
profiles $f_0$ originated at the first critical exponent $p=6^-$;
see Figure \ref{SN1} for the enlarged area of the turning point.}
 \label{FBif2}
\end{figure}

Figure \ref{FBif1}
 illustrates the vanishing behaviour of the VSS profiles $f_0(y)$
as $p \to p_0=n+5$ for $N=1$ and $n=1$. Notice that, according to
(\ref{ep5}), the rate of decay is fast,
 $$
 \|f\|_\infty \sim
(6-p)^{\frac 15} \quad \mbox{as} \quad  p \to 6^-.
 $$
  In Figure \ref{FBif2}, we show the corresponding
first $p$-bifurcation branch that is originated at $p=6^-$. Such a
behaviour of this $p$-branch is similar to that for the semilinear
parabolic equation for $n=0$; cf. \cite{GW2}. The global behaviour
of this branch is  unusual: the branch exhibits a definite
non-monotonicity
and turning for  $p \approx 3.46...\,$.

\begin{figure}
\centering
\includegraphics[scale=0.7]{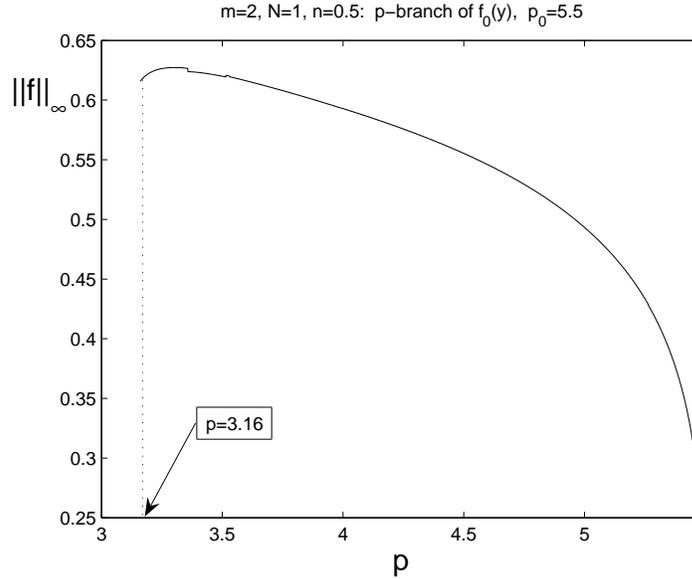}             
 \vskip -.4cm
\caption{\rm\small $N=1$, $n=0.5$: the first $p$-bifurcation
branch  originated at the first critical exponent $p=5.5^-$.}
 \label{FBif3}
\end{figure}

 Such non-monotonicity branching  phenomena were  consistent in numerical experiments. Similar
features are shown for  $n=0.5$ in Figure \ref{FBif3}. Let us
discuss possible consequences  of such a behaviour that was not
available for $n=0$, and hence exists in a strongly quasilinear
case $n>0$ sufficiently large (at least for $n>0.11$ as Figure
\ref{FF12} suggests). It follows from principles of general
branching theory \cite{VaiTr} that if such a $p$-branches
vanishes, it must end up at bifurcation points only.

 On the other hand, we
know that the profiles $f_0$ persist until the critical value
$p=n+1$; see Figure \ref{Def1}, where $f_0(y)$ is available for
$n=0.95$ and $p=2 > n+1=1.95$. Therefore, if a $p_0$-branch
disappear at some bifurcation point $p_k >n+1$, it must appear at
another (possibly, saddle-node) bifurcation point $\hat p_l<p_k$.
As we will explain, at standard pitchfork bifurcations from 0 at
$p=p_l$ with $l=1,2,...$, other types of VSS profiles $f_l$
appear, so that the new appearance of $f_0$
 can be associated with new bifurcations and branching.
Therefore, $f_0$ must appear at some  subcritical saddle-node
bifurcation, most probably embracing branches of profiles $f_0$ and $f_2$
that have a similar geometric structure.

\begin{figure}
 \centering
 \psfrag{||f||}{$\|f\|_\infty$}
 \psfrag{p}{$p$}
 \psfrag{p2}{$p_2$}
 \psfrag{t3}{$t_3$}
 \psfrag{t4}{$t_4$}
  \psfrag{v(x,t-)}{$v(x,T^-)$}
  \psfrag{final-time profile}{final-time profile}
   \psfrag{tapp1}{$t \approx 1^-$}
\psfrag{x}{$x$}
 \psfrag{0<t1<t2<t3<t4}{$0<t_1<t_2<t_3<t_4$}
  \psfrag{0}{$0$}
 \psfrag{l}{$l$}
 \psfrag{-l}{$-l$}
\includegraphics[scale=0.43]{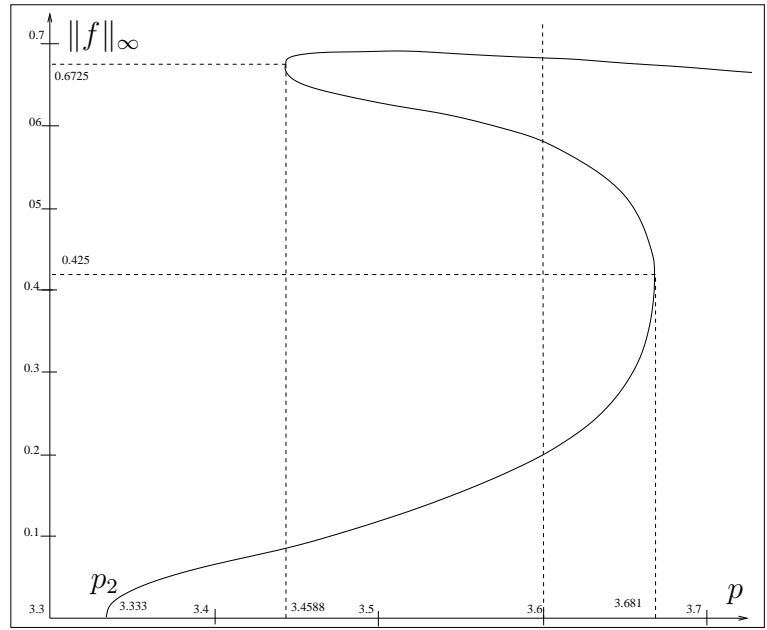}     
\caption{\small The enlarged non-monotonicity  part of the
bifurcation $p_{0,2}$-diagram for the ODE (\ref{OD66}) in $\re$
for $n=1$. The behaviour close to $p_2=3.333...$ with two turning
(saddle-node) points.}
     \vskip -.3cm
 \label{SN1}
\end{figure}
\begin{figure}
\centering
\includegraphics[scale=0.7]{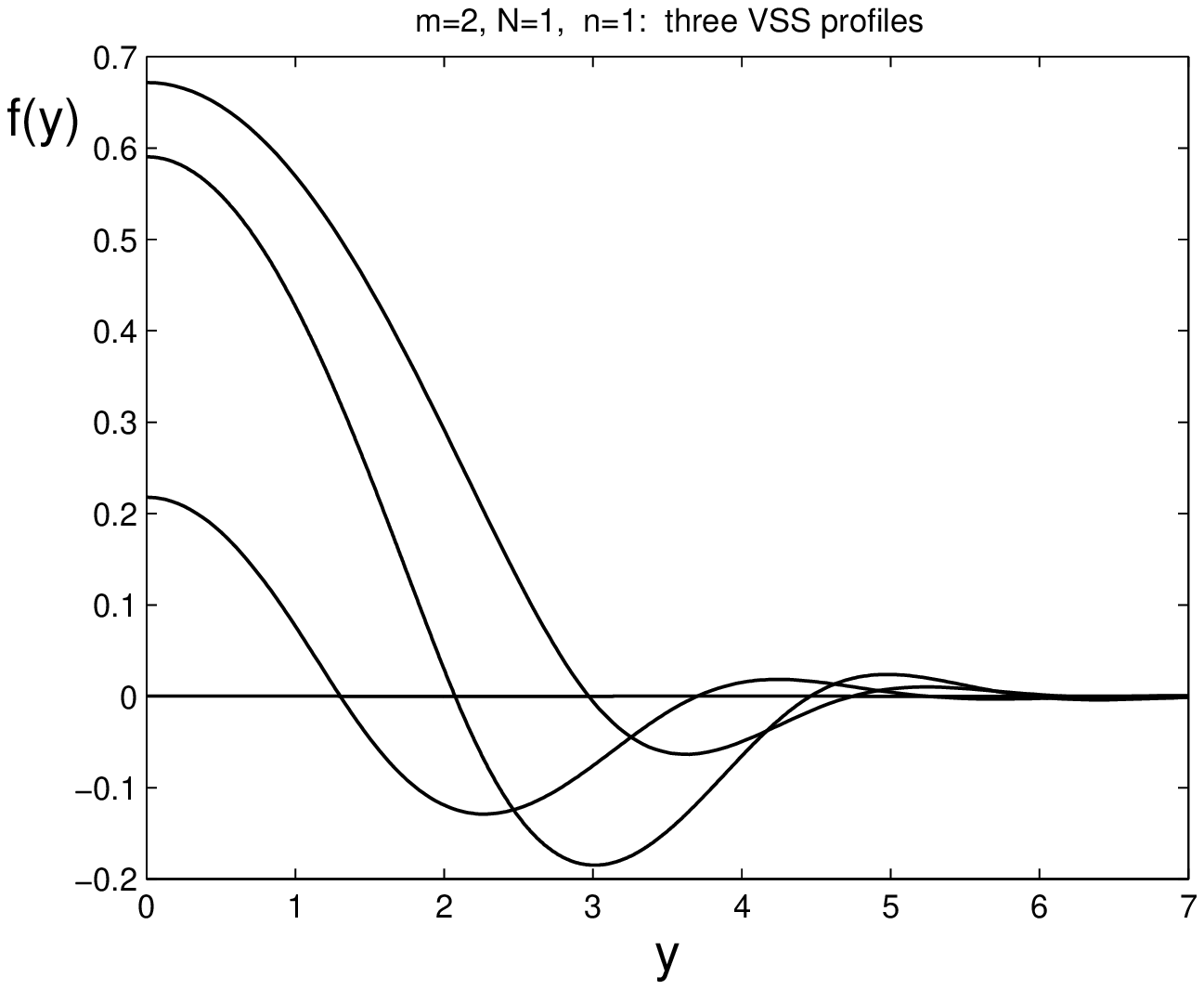}            
 \vskip -.4cm
\caption{\rm\small Three $f_{0,2}(y)$ profiles from Figure
\ref{SN1} for $p=3.6$.}
 \label{FG1}
\end{figure}

We now answer the  question posed above. In Figure \ref{SN1}, we
show the $p_{0,2}$-branch for $n=1$ in a neighbourhood of the
bifurcation point $p_2=3.333...$ and of the turning point shown in
Figure \ref{FBif2}. This branch has two turning points which are
saddle-node bifurcations. It follows that the even profiles
$f_0(y)$ and $f_2(y)$ belong to the same $p_{0,2}$-branch, i.e.,
can be continuously deformed to each other as solutions of the ODE
(\ref{OD66}). Note that the pitchfork bifurcation at $p=p_2$ is
now supercritical (the branch is originated for $p>p_2=3.333...$). Figure
\ref{FG1} shows three different ``$f_0$" profiles for $p=3.6$,
which is shown by the  vertical dashed line in Figure \ref{SN1}.
The smallest profile is actually $f_2(y)$, which thus is homotopic
equivalent to $f_0$ (i.e., admits a homotopic path via a family of operators).

We expect that such  saddle-node bifurcations can occur on other
$p_l$-branches creating necessary profiles in different
$p$-intervals. For instance, we observed an evidence that the next
$p_{1,3}$-branch occurs.

Thus, the $p_{0,2}$-branch is a closed curve.
This type of closed $\mu$-bifurcation branches were earlier found
in \cite[\S~6.4]{GW2} for VSSs with another type of
parameterization in the ODEs like (\ref{OD66}), $n=0$ (so that $\b
= \frac 14$), with the change
 $$
  \mbox{$
  \frac 14 \, f'y \mapsto \mu f' y, \quad \mbox{with parameter} \,\,\, \mu
  \in (0, \frac 14].
 $}
  $$
Then a $\mu$-branch was shown to appear at a pitchfork  bifurcation point and
was continued until another, smaller one, i.e., a  global
continuation of such branches up to $\mu= \frac 14^+$ was
impossible. In addition, essentially non-monotone bifurcation
branches were detected \cite[p.~1802]{BGW1} in rather similar fourth-order ODEs
associated with blow-up in higher-order reaction-diffusion
equations such as
 $$
 u_t = -u_{xxxx} + |u|^{p-1} u \quad \mbox{in} \quad \re \times
 (0,T).
  $$

\smallskip

\noi {\bf On other $p$-bifurcation branches.} A similar
``nonlinear" bifurcation analysis can be performed on the basis of
any suitable similarity profiles $F_l(y)$ of the pure TFE
(\ref{TFE1}). Namely, this profile appears in the similarity
solution of (\ref{TFE1}),
 \be
  \label{ss1}
   \mbox{$
  u_l(x,t)=t^{-\a_l} F_l(y), \quad y= x/t^{\b_l}, \quad \mbox{where} \quad \b_l=
  \frac {1-n\a_l}4>0,
  $}
   \ee
   and $\a_l \in (0, \frac 1n)$ is a parameter. Instead of
   (\ref{Od11}), the function $F_l$ solves the following elliptic
   equation:
 \be
 \label{ss2}
 {\bf B}_l(F) \equiv - \n \cdot (|F|^n \n \D F) + \b_l \n F \cdot
 y + \a_l F=0 \quad \mbox{in} \quad \ren.
  \ee
 According to the principles of self-similarity of the {\em second
 kind} (a notion introduced by Ya.B.~Zel'dovich in 1956
 \cite{Zel56}), the acceptable values of the parameter $\b_l$ are
 chosen from the solvability of the problem (\ref{ss2}), i.e.,
 from
 existence of a nontrivial compactly supported  solution $F_l(y)$ in $\ren$.
  We expect existence of a countable (up to obvious scaling) set of such solutions
  $\Phi=\{\a_l, \, F_l(y)\}$, with most of them not being radially symmetric.
   In a certain sense, this looks like
  a nonlinear extension of the linear eigenvalue problem that
  occurs for $n=0$, where $\b_l= \frac 14$ and  (\ref{ss2}) reads
 \be
 \label{ss3}
  \mbox{$
 {\bf B}_l F \equiv -  \D^2 F + \frac 14 \, \n F \cdot
 y + \a_l F=0 \quad \mbox{in} \quad \ren,
  $}
  \ee
  so that $\l_l \equiv \frac N 4 - \a_l$ are eigenvalues of the non
  self-adjoint operator
   \be
   \label{ss4}
    \mbox{$
    {\bf B} = - \D^2 + \frac 14 \, y \cdot \n + \frac N4 \, I.
    $}
     \ee
This spectrum is discrete and is given by \cite{Eg4}
 \be
 \label{ss5}
  \mbox{$
 \l_l = - \frac l 4, \quad l=0,1,2,... \, ,
  $}
  \ee
and each eigenvalue has finite multiplicity. This determines  all
possible values of the parameters of self-similarity
 \be
 \label{ss6}
 \mbox{$
 \a_l= \frac {N+l}4, \quad l=0,1,2... \, .
  $}
  \ee
For small $n>0$, the linear eigenvalue problem (\ref{ss3}) can
predict the nonlinear eigenfunctions of (\ref{ss2}) by the
$n$-branching approach in the lines similar to that in  Appendix A
in \cite{GHCo}, though the justification in the fourth-order case
is much more difficult.

For larger $n>0$, the nonlinear eigenfunctions $\{F_l\}$ are
unknown even in the ODE case $N=1$ or in radial setting in $\ren$.
Recall that, even in the simplest case $l=1$ meaning the 1-dipole
solution $F_1(y)$, the existence of such a profile for not small
$n$ is still unclear mathematically; see references and results in
\cite{Bow01} and more  recent paper \cite{BW06}.

Finally,  once the nonlinear eigenfunction subset $\{F_l\}$ of the
pure TFE (\ref{TFE1}) is known (say, for small $n>0$), for each
function $F_l(y)$  supported in $B_1$ after rescaling, one can
develop the above bifurcation approach taking in (\ref{ep1})
$\tilde F=F_l$ for any $l=0,1,2,...\,$. This gives the sequence of
critical exponents
 \be
 \label{ss7}
  \mbox{$
   \b \equiv \frac{p-(n+1)}{4(p-1)} = \b_l= \frac {1-n \a_l}4
   \quad
   \Longrightarrow \quad p_l= 1 + \frac 1{\a_l}.
    $}
     \ee
For $n=0$, (\ref{ss7}) and (\ref{ss6}) yield precisely the known
critical bifurcation points (\ref{crN}) studied in \cite{GW2}.

Therefore, a more rigorous bifurcation analysis is available for
$p \approx p_l$ and $n \approx 0$, so this means the
 branching
approach with the parameter $\mu=(p,n)^T \approx (1+ \frac
N{N+l},0)^T$ in $\re^2$.

\subsection{On other VSS profiles}

\begin{figure}
\centering \subfigure[$n= \frac 12, \,p=2$]{
\includegraphics[scale=0.52]{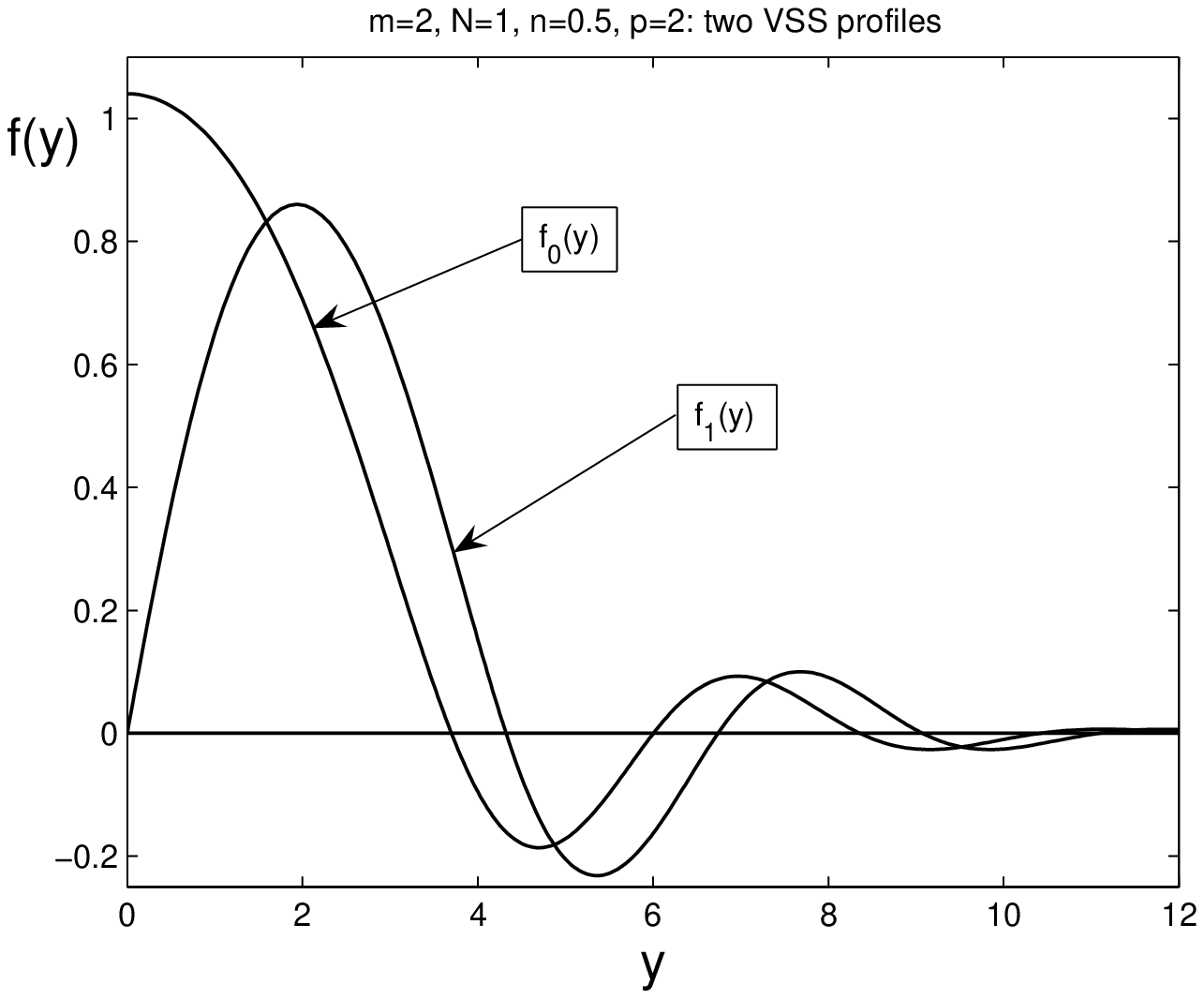}                        
} \subfigure[$n=1, \, p=2.5$]{
\includegraphics[scale=0.52]{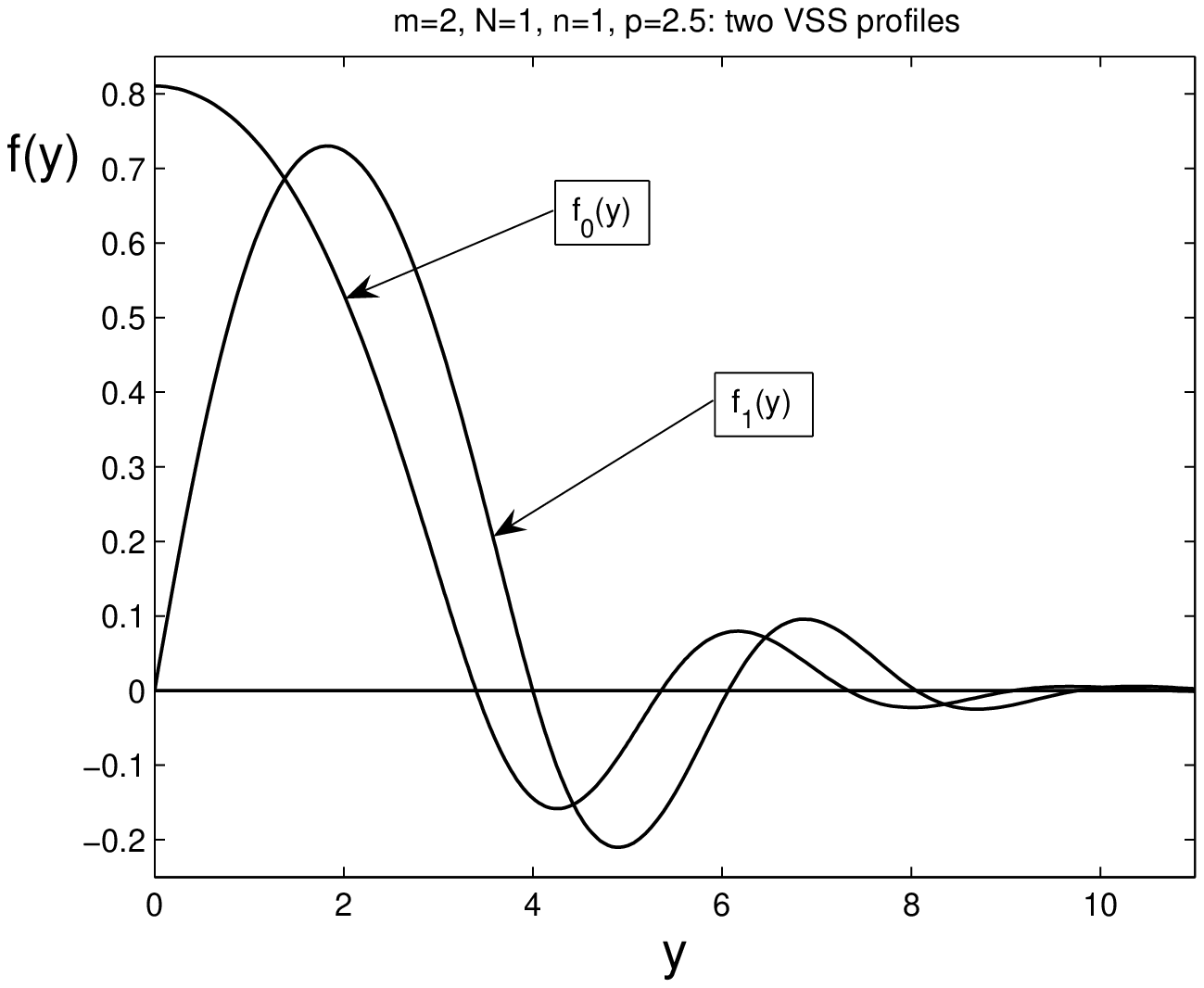}                  
}
 \vskip -.4cm
 \caption{\rm\small First two VSS profiles
  for $n=\frac 12$, $N=1$, $p=2$ (a) and
for $n=1$, $N=1$, $p=2.5$ (b).}
 \label{FBif2V}
\end{figure}

 As we have mentioned, for $n=0$, for any $p$ in the subcritical range $p \in (1,p_0)$, there exists
a finite number of VSS profiles $\{f_l\}$, and each $p$-branch is
originated at bifurcation points (\ref{crN}); see
\cite[\S~6]{GW2}. We claim that similar VSS patterns also exist
for $n>0$ and each $f_{l+1}(y)$ has a more complicated form than
the previous one $f_l(y)$, i.e., has ``more oscillatory" structure
with more non-exponentially small oscillations and ``essential"
zeros. Figure \ref{FBif2V} shows two VSS profiles $f_0$ and $f_1$
for $n=\frac 12$,
 $p=2$ and $n=1$, $p=2.5$ in the 1D case (they look similarly
 emphasizing continuity in $n$).
 It is difficult to check numerically whether other profiles exist.

\begin{figure}
\centering
\includegraphics[scale=0.7]{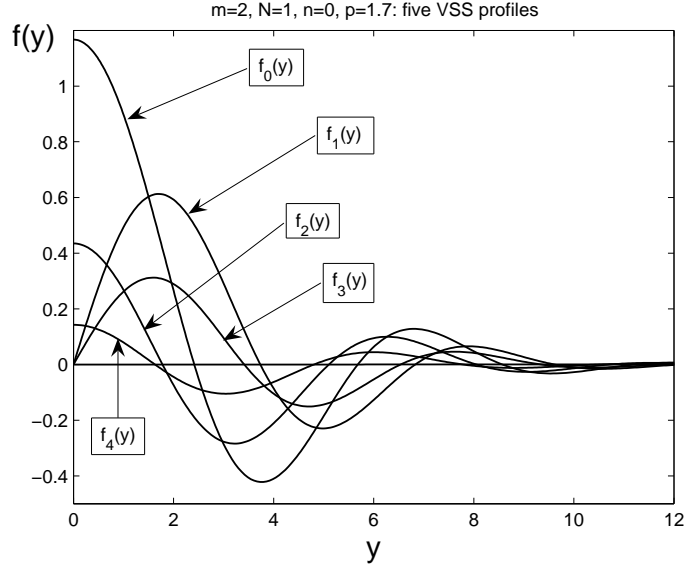} 
 \vskip -.4cm
\caption{\rm\small  Five VSS profiles for the linear diffusion
$n=0$ in 1D, with $p=1.7$.}
 \label{FLin}
\end{figure}

 For comparison, in
 Figure \ref{FLin}, we present five VSS profiles for the
 semilinear PDE for $n=0$, with $N=1$ and $p= 1.7$, which are computed much easier. According to
(\ref{crN}), there exist five bifurcation points above 1.7,
 \be
 \label{ppp1}
  \mbox{$
 p_0=5, \,\,\, p_1=3, \,\,\, p_2=\frac 73, \,\,\, p_3=2, \,\,\,
 p_4= \frac 95= 1.8 \quad (p_5=\frac 53=1.66...<1.7).
  $}
  \ee
Therefore, for $p=1.7$, bearing in mind  that the $p$-bifurcation
branches are monotone decreasing in $p$, there exist precisely
five VSS profiles indicated in this Figure. By continuity, we
expect that these profiles do not essentially change and exist for
all small enough $n>0$, where these become compactly supported
with the interface dependence \cite[\S~10]{Gl4}
 $$
  \mbox{$
 y_0(n) \sim n ^{-\frac 34} \to + \infty \quad \mbox{as} \quad n \to 0^+.
  $}
  $$


\smallskip
{\bf Acknowledgements.}  The author would like to thank P.J.
Harwin for discussions and performing careful  numerical
calculations for Figure 1.


\begin{thebibliography}{10}



   \bibitem
    {Beck05}
  J.~Becker and G.~Gr\"un, {\em The thin-film equation: recent advances and some new perspectives},
  {J.~Phys.: Condens. Matter}, {\bf 17} (2005), S291--S307.


\bibitem 
 {BF1}
 F. Bernis and A. Friedman, \emph{Higher order nonlinear degenerate
 parabolic equations}, J. Differ. Equat. \textbf{83} (1990), 179--206.


 \bibitem 
 {BerHK00}
 F.~Bernis, J.~Hulshof, and J.R.~King,
 {\em Dipoles and similarity solutions of the thin film equation in the half-line},
 {\em  Nonlinearity,} {\bf 13} (2000), 413--439.


  \bibitem 
 {BerHQ00}
 F.~Bernis, J.~Hulshof, and F.~Quir\'os, {\em The ``linear" limit  of thin film flows
as an obstacle-type free boundary problem},
  {SIAM J.~Appl. Math.,} {\bf 61} (2000), 1062--1079.

\bibitem 
 {BMcL91}
 F.~Bernis and J.B.~McLeod, {\em Similarity solutions of a higher order nonlinear
 diffusion equation}, {Nonl. Anal.,} {\bf 17} (1991), 1039--1068.



\bibitem 
 {BPelW92}
 F.~Bernis, L.A.~Peletier, and S.M.~Williams, {\em Source type solutions
 of a fourth order nonlinear degenerate parabolic equation}, {Nonl. Anal.,}
   {\bf 18} (1992), 217--234.






\bibitem 
 {Bow01} M.~Bowen, J.~Hulshof, and J.R.~King, {\em Anomaluous exponents
 and dipole solutions for the thin film equation}, {SIAM J.~Appl.
 Math.,} {\bf 62} (2001), 149--179.

\bibitem 
 {BW06} M.~Bowen and T.P.~Witelski, {\em The linear limit of
 the dipole problem for the thin film equation}, {SIAM J.~Appl.
 Math.,} {\bf 66} (2006), 1727--1748.



\bibitem 
{BGW1}
   C.J.~Budd, V.A.~Galaktionov, and J.F.~Williams, {\em
Self-similar blow-up in higher-order semilinear parabolic
equations}, SIAM J. Appl. Math., {\bf 65} (2004), 1775--1809.

\bibitem
{Car07}
 E.C.~Carlen and S.~Ulusoy, {\em Asymptotic equipartition and long-time
behaviour of solutions of a thin film eqwuation}, J.~Differ.
Equat., {\bf 241} (2007), 279--292.



\bibitem
{CarrT02}
   J.A.~Carrillo and G.~Toscani, {\em Long-time asymptotic
behaviour for strong solutions of the thin film eqwuations}, Comm.
Math. Phys., {\bf 225} (2002), 551--571.




  \bibitem{Eg4}
Yu.V.~Egorov, V.A.~Galaktionov, V.A.~Kondratiev, and
S.I.~Pohozaev,
  \emph{Asymptotic behaviour of global solutions to higher-order semilinear
  parabolic equations in the supercritical range}, Adv. Differ. Equat.,
{\bf 9} (2004), 1009--1038.





\bibitem 
{Ell96} C.M.~Elliott and H.~Garcke, {\em On the Cahn--Hilliard
equation with degenerate mobility}, {SIAM J.~Math. Anal.,} {\bf
27} (1996), 404--423.


\bibitem{EllS}
C.~Elliott and Z.~Songmu, \emph{On the {C}ahn-{H}illiard
equation}, Arch. Rat.
  Mech. Anal., \textbf{96} (1986), 339--357.


   \bibitem{Gl4}
J.D.~Evans, V.A.~Galaktionov, and J.R.~King, 
\emph{Source-type solutions of the fourth-order unstable thin film
equation}, Euro J.~Appl. Math., {\bf 18} (2007), 273--321.


     \bibitem{Bl4}
J.D.~Evans, V.A.~Galaktionov, and J.R.~King, {\em Blow-up
similarity solutions of  the fourth-order unstable thin film
equation},  Euro J.~Appl. Math., {\bf 18} (2007), 195--231.



      \bibitem
      {GBl6}
J.D.~Evans, V.A.~Galaktionov, and J.R.~King, {\em Unstable
sixth-order thin film equation. I. Blow-up similarity solutions;
II. Global similarity patterns},
 {Nonlinearity}, {\bf 20} (2007), 1799--1841, 1843--1881.


\bibitem 
 {BFer97}
 R.~Ferreira and F.~Bernis, {\em Source-type solutions to thin-film equations in higher dimensions},
{European J.~Appl. Math.,} {\bf 8} (1997), 507--524.









  \bibitem
 {GalCr}
  V.A. Galaktionov, {\em Critical global asymptotics in
  higher-order semilinear parabolic equations}, Int.
 J. Math. Math. Sci., {\bf 60} (2003),  3809--3825.




\bibitem{GHCo}
V.A.~Galaktionov and P.J.~Harwin, {\em On evolution completeness
of nonlinear
 eigenfunctions for the porous medium equation in the whole space},
 Advances Differ. Equat., {\bf 10} (2005), 635--674.

 \bibitem
  {PetI}
V.A.~Galaktionov and P.J.~Harwin, {\em On centre subspace
behaviour in thin film equations}, SIAM J. Appl. Math., to appear.


\bibitem
  {AMGV}
 V.A.~Galaktionov and J.L.~V\'azquez, {\rm A Stability Technique  for Evolution Partial Differential Equations.
 A Dynamical Systems Approach}, {\rm Progr. in Nonl. Differ.
 Equat. and their Appl.,} {\bf 56},
Birkh\"auser Boston, Inc., MA, 2004.


\bibitem
 {GSVR} V.A.~Galaktionov and S.R.~Svirshchevskii, {\rm Exact Solutions and
 Invariant Subspaces of Nonlinear Partial Differential Equations in Mechanics and Physics},
  Chapman$\,\&\,$Hall/CRC, Taylor and Francis Group, Boca Raton,
FL,
 2007.



\bibitem{GW2}
V.A.~Galaktionov and J.F.~Williams, \emph{On very singular
similarity solutions of a higher-order
  semilinear parabolic equation}, Nonlinearity, {\bf 17} (2004), 1075--1099.



      \bibitem
      {Gia08}
L.~Giacomelli, H.~Kn\"upfer, and F.~Otto, {\em Smooth
zero-contact-angle solutions to a  thin film equation
 around the steady state},
 {J.~Differ. Equat.,} {\bf 245} (2008), 1454--1506.


\bibitem
{GGK} I.~Gohberg, S.~Goldberg, and M.A.~Kaashoek, {\rm Classes of
Linear Operators}, Vol. {\bf  1}, Operator Theory: Advances and
Applications, Vol. {\bf 49}, Birkh\"auser Verlag, Basel/Berlin,
1990.








\bibitem
 {Govor05}
L.V.~Govor, J.~Parisi, G.H.~Bauer, and G.~Reiter, {\em Instability
and droplet formation in evaporating thin films of a binary
solution}, {Phys. Rev. E}, {\bf 71},  051603
 (2005).


\bibitem 
{Green78} H.P.~Greenspan, {\em On the motion of a small viscous
droplet that wets a surface}, {J.~Fluid Mech.,} {\bf 84} (1978),
125--143.

\bibitem 
{Gr95} G.~Gr\"un, {\em Degenerate parabolic equations of fourth
order and a plasticity model with non-local hardening}, {Z.~Anal.
Anwendungen,} {\bf 14} (1995), 541--573.




\bibitem  
{Ka1}   A.S.~Kalashnikov, {\em Some problems of the qualitative
theory of second-order  nonlinear
 degenerate parabolic equations,}
{Russian Math. Surveys,} {\bf 42} (1987), 169--222.




\bibitem
{KolF} A.N.~Kolmogorov and S.V.~Fomin, {\rm Elements of the Theory
of Functions and Functional Analysis}, Nauka, Moscow, 1976.



\bibitem
   {KrasZ}
M.A.~Krasnosel'skii and P.P.~Zabreiko, {Geometrical Methods of
Nonlinear
  Analysis}, Springer-Verlag, Berlin/Tokyo, 1984.




\bibitem
{LPugh} R.S. Laugesen and M.C. Pugh, {\em Energy levels of steady
states for thin-film-type equations}, J. Differ. Equat., {\bf 182}
(2002), 377--415.




 \bibitem 
 {Nai1}
 M.A.~Naimark, {\rm Linear Differential Operators}, Part II, Frederick Ungar
Publ. Co.,
 New York, 1968.



\bibitem 
 {Oron97} A.~Oron, S.H.~Davies, and S.G.~Bankoff, {\em Long-scale
 evolution of thin liquids films}, {Rev. Modern Phys.,}
  {\bf 69} (1997), 931--980.






\bibitem 
{Smyth88} N.F.~Smyth and J.M.~Hill, {\em High-order nonlinear
diffusion}, {IMA J.~Appl. Math.,} {\bf 40} (1988), 73--86.

\bibitem 
{VaiTr} M.A. Vainberg and V.A. Trenogin, {\rm Theory of Branching
of Solutions of Non-Linear Equations}, Noordhoff Int. Publ.,
Leiden, 1974.




\bibitem{WitBerBer}
T.P.~Witelski, A.J.~Bernoff, and A.L.~Bertozzi, \emph{Blow-up and
dissipation
  in a critical-case unstable thin film equation},  Euro J.~Appl.
  Math., {\bf 15} (2004), 223--256.



\bibitem 
{Zel56} Ya.B.~Zel'dovich, {\em The motion of a gas under the
action of a short term pressure shock}, Akust. Zh., {\bf 2}
(1956), 28-38; Soviet Phys. Acoustics, {\bf 2} (1956), 25--35.




\end{thebibliography}

\end{document}